\documentclass[11pt]{amsart}
\usepackage{amsmath,amssymb,amsthm}
\usepackage{enumerate}
\usepackage[english]{babel}
\usepackage{pstricks}
\usepackage{graphicx}
\usepackage{curve2e}
\usepackage{tikz}
\usetikzlibrary{matrix}

\usepackage[pagewise]{lineno}

\makeatletter
\@namedef{subjclassname@2010}{%
  \textup{2010} Mathematics Subject Classification}
\makeatother

\theoremstyle{plain}
\newtheorem{main-theorem}{Theorem}

\newtheorem{Thm}{Theorem}

\newtheorem{lemm}{Lemma}[section]
\newtheorem{prop}[lemm]{Proposition}
\newtheorem{coro}[lemm]{Corollary}
\theoremstyle{definition}
\newtheorem{defi}[lemm]{Definition}

\newtheorem{rema}[lemm]{Remark}

\newtheorem*{clai-nn}{Claim}
\newtheorem*{que*}{Major Quesion}
\newtheorem*{rem*}{Remark}
\newtheorem*{theorem*}{Theorem}
\newtheorem*{lem*}{Lemma}

\newcommand{\bbZ}{\mathbb{Z}}
\newcommand{\bbQ}{\mathbb{Q}}
\newcommand{\bbR}{\mathbb{R}}
\newcommand{\bbS}{\mathbb{S}}
\newcommand{\bbT}{\mathbb{T}}
\newcommand{\bbC}{\mathbb{C}}
\newcommand{\bbD}{\mathbb{D}}

\newcommand{\Cc}{\mathcal{C}}

\newcommand{\Lc}{\mathcal{L}}
\newcommand{\Ic}{\mathcal{I}}

\newcommand{\Hc}{\mathcal{H}}

\newcommand{\Dc}{\mathcal{D}}
\newcommand{\M}{\mathcal{M}}

\newcommand{\Oc}{\mathcal{O}}
\newcommand{\Pc}{\mathcal{P}}

\newcommand{\Rc}{\mathcal{R}}
\newcommand{\Tc}{\mathcal{T}}

\newcommand{\Wc}{\mathcal{W}}

\newcommand{\Uc}{\mathcal{U}}
\newcommand{\pie}{\widetilde{\pi}}
\makeatletter
\@namedef{subjclassname@2020}{\textup{2020} Mathematics Subject Classification}
\makeatother

\begin{document}
\title[]{To Define the Core Entropy for  All Polynomials Having a Connected Julia Set}

\author[J. Luo]{Jun Luo} \address{School of Mathematics\\     Sun Yat-Sen University\\ Guangzhou 510275, China} \email{luojun3@mail.sysu.edu.cn}

\author[B. Tan]{Bo Tan} \address{School of Mathematics and Statistics, Huazhong University of Science and Technology, Wuhan 430074, China} \email{tanbo@hust.edu.cn}

\author[Y. Yang]{Yi Yang} \address{School of Mathematics(Zhuhai), Sun Yat-sen University, Zhuhai, 519082 China} \email{yangy699@mail.sysu.edu.cn\ (corresponding author)}

\author[X.T. Yao]{Xiao-Ting Yao} \address{School of Mathematics and Statistics, Guangdong University of Technology, Guangzhou 510520, China} \email{yaoxiaoting55@gdut.edu.cn}

\date{}

{\footnotesize
\begin{abstract}
For all polynomials $f$ with
${\rm deg}(f)\ge2$  that have a connected filled Julia set $K$, we introduce a new quantity $h_{\rm GCE}(f)$, such that  $h_{\rm GCE}\left(f^n\right)=n\cdot h_{\rm GCE}(f)$ for all $n\ge1$ and $h_{\rm GCE}(f)=h_{\rm GCE}(g)$ for $J$-equivalent  $f$ and $g$.  When the coefficients and the critical points of $f$ are real, $h_{\rm GCE}(f)=h(K\cap\bbR,f)$. When $f$ is post-critically finite,  $h_{\rm GCE}(f)$ equals the core entropy $h(\Hc(f),f)$,  where  $\Hc(f)$ is the Hubbard tree. For $f_c(z)=z^2+c$ with $c$ varying in the Mandelbrot set $\M$, the entropy map $c\mapsto h_{\rm GCE}(f_c)$ is not continuous. However, its lower envelope
$h_{\rm core}:\M\rightarrow\bbR$ given by $h_{\rm core}(c)=\inf\left\{t:\ \exists\ c_n\ne c\  \text{with}\ c_n\rightarrow c\ \text{and}\ t=\lim\limits_{n\rightarrow\infty}h_{\rm GCE}\left(f_{c_n}\right)\right\}$ is continuous over $\M$ and has three properties. First, every $h_{\rm core}^{-1}([0,s])$ with $s\ge0$  is connected. In particular, $h_{\rm core}^{-1}(0)$ coincides with  the central molecule of $\M$. Second, $h_{\rm core}(c)=h(\bbR,f_c)$  for $c\in[-2,\frac14]$. Third, $h_{\rm core}(c)=h(\Hc(f_c),f_c)$  for post-critically finite $f_c$.
\end{abstract}
}


\maketitle

\vspace{-0.382cm}

{\footnotesize
\tableofcontents
}

\newpage

\section{Introduction and Main Results}

We discuss polynomials $f$  with degree $d\ge2$, to be considered as self-maps of the extended complex plane $\widehat{\bbC}$. Every polynomial is conjugated by an affine map to a  monic centralized one, of the form $z^d+\sum_{i=0}^{d-2}a_iz^i$. So, we may as well specialize to  such polynomials. Our aim is to introduce and study the {\bf generalized core entropy} (shortly, GCE). This new quantity, denoted by $h_{\rm GCE}(f)$, is universally defined for all $f$ having a connected filled Julia set $K$.  Here we do not require that $K$ be locally connected and even allow $f$ to have  Crem\'er points or Siegel disks. 

We will illustrate that  $f\mapsto h_{\rm GCE}(f)$ fully generalizes two similar maps from the literature that are related to the entropy of $f$ restricted to certain invariant subsets. 
First, $h_{\rm GCE}(f)$ equals the core entropy $h(\Hc(f),f)$, whenever $f$ is a post critically finite (shortly, PCF) polynomial. 
Second, $h_{\rm GCE}(f)=h(\bbR,f)$, if only the  critical points of $f$ are all real. 
We also identify  a subset $B(f)$ of $\bbS^1=\{z\in\bbC: |z|=1\}$ whose Hausdorff dimension, under mild assumptions, coincides with the ratio $\frac{h_{\rm GCE}(f)}{\log d}$. 
For quadratics $f_c(z)=z^2+c$, the map $c\mapsto h_{\rm GCE}(f_c)$ as a function  of the Mandelbrot set $\M$  is not continuous, although it is continuous at the points of the interior $\M^o$ and  at $c=-2$. However, the lower envelope of $c\mapsto h_{\rm GCE}(f_c)$ is naturally defined and continuous on the whole Mandelbrot set. Moreover, it extends the classical core entropy map $c\mapsto h(\Hc(f_c),f_c)$, defined on  $\{c\in\M:f_c\ \text{is\ PCF}\}$.

Our motivations come from very recent studies, especially \cite{Bruin-vanStrien15,Dudko-Schleicher20,Gao20,Gao-Tiozzo22,Jung14,Milnor-Thurston,Milnor-Tresser00,Tiozzo15,Tiozzo16}. Before diving into further details, let us  recall notions and known results from the literature.

Due to \cite{Lyubich83}, the topological entropy of $f$, as a self-map of $\widehat{\bbC}$, is known to be constantly equal to $\log d$. Since $f$ restricted to its Fatou set is equicontinuous, the entropy of $f$  restricted to its Julia set $J$, denoted by  $h(J,f)$, also equals $\log d$.  Therefore, no two polynomials of the same degree can be told apart by checking the  entropy $h(J,f)$ alone. 

If $f$ is post-critically finite (shortly, PCF), its filled Julia set $K$ is connected and locally connected. And one may restrict $f$ to the Hubbard tree $\Hc(f)$, a minimal tree in $K$ that is  invariant under $f$ and contains the critical points of $f$. Here the edges of $\Hc(f)$ are required to be regulated arcs, so that the choice of $\Hc(f)$ is actually unique \cite[Lemma 2.1]{Poirier10}. In the literature, the Hubbard tree sometimes means the subtree $f(\Hc(f))\subset\Hc(f)$. Since every invariant probability measure of $(\Hc(f),f)$ is supported on $f(\Hc(f))$, it makes no difference from the viewpoint of ergodic theory if one restricts instead to the subsystem $\Big(f(\Hc(f)),f\Big)$.  

 Due to \cite{BFH92,Poirier10}, it is known that  PCF polynomials  are $J$-equivalent if and only if they share the same critical portrait.  Therefore, the subsystem $(\Hc(f),f)$ and its entropy $h(\Hc(f),f)$ are very useful, when one wants to tell apart PCF polynomials of the same degree. In deed, the analysis of $h(\Hc(f),f)$ has been a major theme in recent studies such as \cite{Gao20,Gao-Tiozzo22,Tiozzo15,Tiozzo16}.
By  \cite{Tiozzo16} one can find a continuous function $h_T:\M\rightarrow\bbR$, satisfying $h_T(c)=h\left(\Hc(f_c),f_c\right)$ for all $c\in\M$ such that $f_c(z)=z^2+c$ is PCF. The  function $h_T$ relies on a continuous map $\theta\mapsto\log r_\theta$ with $\theta\in\mathbb{R}/\mathbb{Z}$, given in  \cite[Theorem 8.2]{Tiozzo16} and independently in \cite[Theorem 1.1]{Dudko-Schleicher20}. The existence of such a function $h_T:\M\rightarrow\bbR$ was conjectured by Thurston and known to specialists. For the sake of convenience, we will refer to it as {\bf Thurston's core entropy map}. 

Two further properties of $h_T:\M\rightarrow\bbR$  are noteworthy. First, $h_T(c)=h(\bbR,f_c)$ for $c\in[-2,\frac14]$. Second, $h_T(c_1)\le h_T(c_2)$ for  $c_1,c_2\in \Dc$ that satisfy $c_1<c_2$ in the sense of  \cite[$\S20.4$]{DH84}, where $\Dc$ consists of the roots and the centers of hyperbolic components and Misurewicz points. Since  $\overline{\Dc}\supset\partial\M$,  $h_T$ continuously extends the following well-known entropy map. 
\begin{theorem*}[\bf{\cite[Lemma 12.3 and Corollary 13.1]{Milnor-Thurston}}]
 Let $\displaystyle I_c=\left[\frac12-\sqrt{\frac14-c},\frac12+\sqrt{\frac14-c}\right]$ and $f_c(z)=z^2+c$  for $c\in[-2,\frac14]$. Then $c\mapsto h(I_c,f_c)$ is a nonincreasing continuous function.
 \end{theorem*}

 As was done in earlier works, such as \cite{Tiozzo15},  one may analyze $\M$ by investigating the level sets of $h_T$ and their components. Notice that, for $c\in\M$ such that $f_c(z)$ is not PCF, the value $h_T(c)$ is not directly connected to the dynamics of $f_c$. In deed,  if the filled Julia set $K_c$ is not locally connected (even not path connected), no invariant subtree of $K_c$ containing the critical point has been found. Such an issue was already observed in \cite[Remark 4.5]{Jung14}. 
Therefore, one needs to figure out whether and how  $h_T(c)$  is connected to the dynamics of $f_c$ when  $K_c$ is  not  locally connected, especially when there is a Crem\'er point or a Siegel disk. 

The map $f\mapsto h_{\rm GCE}(f)$ we study will provide a reasonable resolution to the previous issue. The quantity $h_{\rm GCE}(f)$ will be universally defined for any $f$ that has a connected Julia set. Focusing on quadratics $f_c(z)=z^2+c$ , we will analyze the map $c\mapsto h_{\rm GCE}(f_c)$ with $c\in\M$. This is a bounded function whose discontinuities form a dense subset of $\partial\M$. However, its lower envelope $h_{\rm core}:\M\rightarrow\bbR$, to be defined later on, satisfies $h_{\rm core}(c)=h_T(c)$ for all $c\in\M$. Therefore, we can  directly identify $h_T$, in a way that is  independent of the two approaches given in  \cite{Dudko-Schleicher20,Tiozzo16}. When proving the continuity of $h_{\rm core}$, we still need to use the continuity of the map $\theta\mapsto\log r_\theta$  given in \cite[Theorem 8.2]{Tiozzo16}. Here one may wonder how to find a self-contained proof for the continuity of $h_{\rm core}$, without using that of $\theta\mapsto\log r_\theta$.


The notion of $h_{\rm GCE}(f)$ benefits from earlier works that identify and study  special sub-systems of $(K,f)$, when $f$ PCF and when its critical points are all real.
In the first,  one may  restrict $f$ to  the Hubbard tree $\Hc(f)$ and consider the subsystem $(\Hc(f),f)$. 
In the second,  $K\cap\bbR$ is  an interval and the entropy map $f\mapsto h(K\cap\bbR,f)$  demonstrates very rich properties \cite{Bruin-vanStrien15, Douady95, Milnor-Thurston,Milnor-Tresser00}.

\begin{defi}\label{def:real_classical}
When $f$ is PCF, we call  $h(\Hc(f),f)$  the {\bf classical core entropy} of $f$. When the critical points of $f$ are all real, we call $h(K\cap\bbR,f)$  the {\bf real entropy} of $f$.
\end{defi}

Notice that if $f$ is a polynomial whose critical points are all real then  the minimal interval $I_f\subset K\cap\bbR$ containing the critical points satisfies $f(I_f)\subset I_f$. If in addition $f$ is PCF, we  have $I_f=\Hc(f)$. Thus $h(\Hc(f),f)$ is a variant of $h(\bbR,f)$ for  polynomials whose critical points are all real and applies to complex polynomials that are PCF. 

In the current paper, we aim to extend the notion of $h(\Hc(f),f)$ so that it may fit well with non-PCF polynomials and hence generalize the real entropy in the same time. To do that,  one may remove the PCF assumption and require instead that there be a minimal invariant tree $T_f\subset K$ containing all the critical points. The existence of such a tree $T_f$ usually requires  $K$ to be locally connected (at least path connected), which may not be the case. To bypass this seemingly obstacle, we turn to consider {\bf dendrite factors} $(T,g)$ of $(J,f)$, where $g$ is a continuous self-map of a dendrite and the underlying semi-conjugacy $\phi:J\rightarrow T$ is a monotone map, in the sense that the point inverses $\phi^{-1}(u)$ are each a sub-continuum of $J$. 
\begin{defi}\label{def:MDF}
A  dendrite factor  $\left(T_0,g_0\right)$ of $(J,f)$ is said to be {\bf maximal} provided that every dendrite factor of $(J,f)$ is also a dendrite factor of $\left(T_0,g_0\right)$.
\end{defi}

The first main result is the following.
\begin{Thm}\label{theo:MDF} 
$(J,f)$ has a maximal dendrite factor (shortly, MDF). \end{Thm}

The MDF is unique. In deed, every  dendrite factor gives rise to an upper semi-continuous (shortly, USC) decomposition of $J$ into sub-continua, which are exactly the point inverses of the underlying semi-conjugacy. Theorem \ref{theo:C_for_core} will show that among the family of such decompositions there is a finest one, which is determined by the topology of $J$.  This neatly resolves the existence and the uniqueness of the MDF. Hereafter, we denote by $\left(\Tc(f),\overline{f}\right)$ the MDF of $(J,f)$ and by $\pie:J\rightarrow\Tc(f)$ the underlying semi-conjugacy.

Theorem \ref{theo:dichotomy} will show that either $\Tc(f)$ is a point or  $h\left(\Tc(f),\overline{f}\right)=\log d$.  In Lemmas \ref{lem:degenerate_core} and \ref{lem:degenerate_core_converse} we restrict to quadratics $f_c(z)=z^2+c$ with $c\in\M$ and prove that $\Tc(f_c)$ is a point if and only if there are  hyperbolic components $H_1,\ldots,H_k$ whose closures form a continuum containing $\{0,c\}$.  Such a parameter $c$ lies in the {\bf  central molecule} $M_0$ \cite{KahnLyubich09-a}. 
By Corollary \ref{cor:tuning_molecule-a} we know that there are infinitely many $c'\in M_0$ (including the Feigenbaum parameter $c_F$)  such that   $\Tc_0\left(f_{c'}\right)$ is not a point and   $h_{\rm GCE}\left(f_{c'}\right)=0$. 


To identify $h_{\rm GCE}(f)$, we may pick special sub-systems of $\left(\Tc(f),\overline{f}\right)$ in the following way. 
Let  ${\rm Crit}(f)$ be the set of the critical points and $P_f$ the {\bf post critical set}, which is the closure of $\Big\{f^n(x): x\in{\rm Crit}(f), n\ge1\Big\}$. Then there are minimal sub-continua  $\Tc_0(f),\Tc_1(f)\subset\Tc(f)$  that are invariant under $\overline{f}$ and respectively contain $\pie\left({\rm Crit}(f)\right)$ and $\pie\left(P_f\right)$. 
The resulting sub-system $\left(\Tc_0(f),\overline{f}\right)$ is called  the {\bf dynamic core} of $f$ and  $\left(\Tc_1(f),\overline{f}\right)$  the {\bf reduced dynamic core}. We also call  $\Tc_0(f)$ the {\bf topological core} and $\Tc_1(f)$ the {\bf reduced topological core}.  It is immediate that $h\left(\Tc_0(f),\overline{f}\right)=h\left(\Tc_1(f),\overline{f}\right)$, since $\overline{f}\left(\Tc_0(f)\right)=\Tc_1(f)$. 
\begin{defi}\label{def:core}
Set $h_{\rm GCE}(f)=h\left(\Tc_0(f),\overline{f}\right)$ and call it the {\bf  generalized  core entropy} of $f$.
\end{defi}

There is another subsystem, which is useful  when $f$ is renormalizable. Let $M_f$ be the minimal compact set such that $\mu(M_f)=1$ for each invariant probability $\mu$ of $(\Tc_0(f),\overline{f})$. Let $\Tc_{\rm irr}(f)\subset\Tc_0(f)$ be the minimal sub-continuum containing $M_f$ and  call it the {\bf irreducible topological core}. Clearly, $\Tc_{\rm irr}(f)$ is invariant under $\overline{f}$. 
\begin{defi}\label{def:irreducible_core}
We call $\left(\Tc_{\rm irr}(f),\overline{f}\right)$  the {\bf irreducible dynamic core} of $f$. \end{defi}


For each $d\ge2$, the map $f\mapsto h_{\rm GCE}(f)$ is  a well-defined map of the  connectedness locus $\Cc_d$ that consists of the coefficient vectors $\textbf{a}=(a_0,\ldots,a_{d-2})\in\bbC^{d-1}$ such that
$f_\textbf{a}(z)=z^d+\sum_{j=0}^{d-2}a_jz^j$ has a connected Julia set. The connectedness locus $\Cc_2$ is well known as the Mandelbrot set and often denoted by $\M$. Due to \cite{DH82}, $\M$ is a cellular continuum. A long standing open question is the local connectedness of $\M$. By \cite[Proposition 3.1]{Branner93}, every $\Cc_d(d\ge3)$ is also a cellular continuum.
 By \cite[Theprem 3.2]{Branner93}, $\Cc_d(d\ge3)$ is  NOT locally connected.

For the sake of convenience, we often identify coefficient vectors $\textbf{a}\in\Cc_d$ with polynomials $f_\textbf{a}(z)$ and vice versa.
If $\textbf{a}\in\bbR^d$,  the filled Julia set of $f_\textbf{a}$ intersects $\bbR$ at an interval $[a,b]$ such that $f(\{a,b\})\subset\{a,b\}$. Thus $\left.f\right|_{[a,b]}$ is a boundary anchored $(d-1)$-modal map \cite[\S2]{Milnor-Tresser00}.

The second main result is the following.


\begin{Thm}\label{theo:good_extension}
If ${\rm Crit}(f)\subset\bbR$, $h_{\rm GCE}(f)=h(K\cap\bbR,f)$. If $f$ is PCF, $h_{\rm GCE}(f)=h(\Hc(f),f)$.
\end{Thm}

\begin{rema}
The new quantity $h_{\rm GCE}(f)$ fully generalizes both the real entropy and the classical core entropy. There are some basic observations.
First,  if a polynomial-like map $f$ with Julia set $J$ is hybrid equivalent to a polynomial $f_0$ having a connected Julia set $f_0$, then  $(J,f)$ is topologically conjugate with  $(J_0,f_0)$. It follows that the MDF's of $(J,f)$ and $(J_0,f_0)$ are topologically conjugate, too. Therefore, 
$h_{\rm GCE}(f)$ is well-defined and coincides with $h_{\rm GCE}(f_0)$.  Second, by definition $h_{\rm GCE}(f)=h_{\rm GCE}(g)$ is true  for all  $g$ that are $J$-equivalent to $f$. Third, it is routine to check that $h_{\rm GCE}(f^n)=n\cdot h_{\rm GCE}(f)(n\ge1)$.
Finally, if a polynomial-like map  $f$ has a connected filled  Julia set $K$ and $f^p(p\ge2)$ restricted to a neighborhood $U$ of a continuum $Y\subset K$ is a polynomial-like map hybrid equivalent to a polynomial $f_0$ restricted to a neighborhood $U_0$ of its filled  Julia set $K_0$, then the irreducible dynamic core of $f_0$ is topologically conjugate with a sub-system of the dynamic core of $f^p$. It follows that $h_{\rm GCE}(f)\ge\frac1p h_{\rm GCE}(f_0)$.
\end{rema}


The last main result concerns quadratics $f_c(z)=z^2+c$ and  the map  $c\mapsto h_{\rm GCE}(f_c)$  with  $c\in\M$.  Let us introduce further terminology below. 

\begin{defi}\label{def:cluster}
The {\bf cluster set} of $h_{\rm GCE}:\Cc_d\rightarrow\bbR$ at $f\in\Cc_d$, denoted by  $C(h_{\rm GCE},f)$,  is  the set of  $t\in\bbR$ such that there are $f_n\in\Cc_d\setminus\{f\}$ satisfying $\lim\limits_{n\rightarrow\infty}h_{\rm GCE}(f_n)=t$ and $f_n\rightarrow f$, in the sense that the coefficient vectors of $f_n$ converge to that of $f$. 
\end{defi}

\begin{defi}\label{def:envelope}
Set $h_{\rm core}(f)=\inf C(h_{\rm GCE},f)$ and call it the {\bf lower envelope} of $h_{\rm GCE}:\Cc_d\rightarrow\bbR$. \end{defi}

Let us mention some basic facts concerning the case $\Cc_2=\M$. First, by \cite{MSS83} the map $h_{\rm GCE}:\M\rightarrow\bbR$ is constant on every component of the interior $\M^o$ hence is continuous at each point of $\M^o$.  Second,  $h_{\rm GCE}^{-1}(\log2)$ is a dense subset of the boundary $\partial\M$. Here we note that for $c$ belonging to the impression ${\rm Imp}(\M,\theta)$ such that 
$\Big\{2^n\theta({\rm mod}1)\Big\}$ is dense in $\bbR/\bbZ$,  the topological core $\Tc_0(f_c)$ is just the dendrite $\Tc(f_c)$ hence $h_{\rm GCE}(f_c)=\log2$. It follows that the previous  dense subset of $\partial\M$ has full harmonic measure and that $h_{\rm GCE}:\M\rightarrow\bbR$ is  continuous at none of the parameters $c\in\partial\M$ with $h_{\rm GCE}(f_c)<\log2$. Such parameters $c$ include all those that lie on the boundary of a hyperbolic component. Therefore, they form a dense subset of $\partial\M$.  Third, $h_{\rm GCE}:\M\rightarrow\bbR$ is continuous at $c=-2$ by \cite[Lemma 12.3]{Milnor-Thurston},  Theorem \ref{theo:good_extension} and  Theorem \ref{theo:monotonicity}.  

However, the lower envelope $h_{\rm core}:\M\rightarrow\bbR$ is continuous over the whole Mandelbrot set $\M$.
\begin{Thm}\label{theo:properties_of_h}
$h_{\rm core}^{-1}([0,t])$ is a continuum for all $t\ge0$ and $h_{\rm core}\!:\M\rightarrow\bbR$ is continuous. \end{Thm}

\begin{rema} Note that, by \cite[Theorem 1.2]{Gao-Tiozzo22}, the classical core entropy map in any degree $d\ge3$   depends continuously on the coefficient vector.  One may wonder whether there is a continuous extension over the whole parameter space $\Cc_d$ (or certain  subspaces). Also note that our proof for the continuity of $h_{\rm core}:\M\rightarrow\bbR$ depends on   \cite[Theorem 8.2]{Tiozzo16} and uses some well known topological properties of $\M$. The same proof may not fit well with higher degree polynomials. However, if one focuses on $\M_d(d\ge3)$, which consists of all  $\textbf{a}\in\Cc_d$ with $a_1=\cdots=a_{d-2}=0$, a continuous extension of the classical core entropy map may be obtained in a similar way. Discussions along such directions can be expected.  \end{rema}


Theorems \ref{theo:MDF} to \ref{theo:properties_of_h} are the main results. Section \ref{sec:general} prepares the general theory of GCE and proves Theorems \ref{theo:MDF} and \ref{theo:good_extension}. Section \ref{sec:quadratic} focuses on quadratics, analyzes the GCE map $c\mapsto h_{\rm GCE}(f_c)$ and proves Theorem \ref{theo:properties_of_h}.  In sections \ref{sec:general} and \ref{sec:quadratic} we also obtain some results, Theorems \ref{theo:invariance_sim} to \ref{theo:continuity_h_core}, that are  fundamental in our discussions. All the other results  appear as lemmas, propositions, or corollaries. Together with definitions and remarks, they will be numbered within sections. 

\section{General Theory of Generalized Core Entropy (GCE)}\label{sec:general}
The materials of this section are outlined as follows. Subsection \ref{section:MDF} gives a complete proof for Theorem \ref{theo:MDF}. Subsection \ref{section:doubling_dynamics}  addresses how the dynamics of $\left(\Tc(f),\overline{f}\right)$ is related to that of $\left(\bbS^1,\sigma_d\right)$, where  $\sigma_d(z)=z^d$. Among others, a dichotomy of $\left(\Tc(f),\overline{f}\right)$ is obtained in Theorem \ref{theo:dichotomy}.   Subsection \ref{section:biac_dim}  identifies  a special set $B(f)\subset\bbS^1$  for  any polynomial $f\in\Cc_d \ (d\ge2)$ whose Hausdorff dimension, under mild conditions,  equals the ratio $\displaystyle\frac{h_{\rm GCE}(f)}{\log d}$. See Theorem  \ref{theo:biac_dim}. Subsection \ref{section:good_extension} proves Theorem \ref{theo:good_extension}. Subsection  \ref{section:small_J}  deals with renormalizable polynomials and illustrates how $h_{\rm GCE}(f)$ is related to the dynamics of the ``small Julia set''.

\subsection{The Existence of the Maximal Dendrite Factor (MDF)} \label{section:MDF}
We  shall need  the core decomposition $\Dc_{K}^{PC}$ with  Peano quotient for  nonempty planar compacta $K\subset\widehat{\bbC}$.  This  is the finest USC decomposition of $K$ into sub-continua such that the quotient space is a Peano compactum.
Here  a {\bf Peano compactum} is a compact metric space with locally connected components such that for any $C>0$ at most finitely many of its components are of diameter $>C$.
See \cite[Theorem 7]{LLY19} for the existence of such a core decomposition. See  \cite[Theorem 3]{LLY19} for the motivation of Peano compactum from the  two dynamic models discussed in \cite{BCO11,BCO13}, which are helpful in the analysis of polynomial Julia sets.



Recall that $\Dc_K^{PC}$ is an {\bf atomic decomposition} of $K$ in the sense of \cite[Definition 1.1]{FitzGerald67}, with respect to the property of having a Peano compactum as the quotient space.
Thus elements of $\Dc_{K}^{PC}$  are also called {\bf atoms} of $K$.
By  \cite[Theorem 1.1]{LYY25}, for any compact set $K\subset\widehat{\bbC}$ and any  rational map $f$,   every atom of $f^{-1}(K)$ is sent by $f$ onto an atom of $K$.
Here $K$ may not be the Julia set of $f$.  When $f$ is a polynomial and $K$ the filled  Julia set of $f$, such an invariance is already known from  \cite{BCO11,BCO13}. We will see how this invariance plays a role in  the theory of GCE.

\begin{defi}\label{def:sim_f}
Given a polynomial $f$ with filled  Julia set $K$, let $\Dc(f)$ be the finest USC decomposition that splits no atom of $K$  and none of the bounded Fatou components; moreover, let $\sim_f$ be the equivalence on $\widehat{\bbC}$ whose classes are exactly the elements of $\Dc(f)$.
\end{defi}

Every $\delta\in\Dc(f)$ is either a singleton disjoint from  $K$ or a continuum (possibly a point) which is contained in $K$ and has a connected complement in $\widehat{\bbC}$. By Moore's Decomposition Theorem \cite[Theorem 22]{Moore25}, the resulting quotient space of $\Dc(f)$ is homeomorphic with $\widehat{\bbC}$. Thus the projection   $\pie: \widehat{\bbC}\rightarrow\Dc(f)$ is a branched cover and $\Tc(f)=\pie(K)$ a planar Peano compactum. Since $\pie\left(\widehat{\bbC}\setminus K\right)$ is open and dense in the quotient space, every component of  $\Tc(f)$ is a dendrite. 
\begin{main-theorem}\label{theo:invariance_sim}
For each $\delta\in\Dc(f)$ we have two assertions: (1)  $f(\delta)\in\Dc(f)$; (2) $\delta\cap J$ is connected.\end{main-theorem}

\begin{proof}
Let $R_f$ be the relation on $\widehat{\mathbb{C}}$ that identifies $x,y$  which are are contained in a single atom of $K$  in the closure of a bounded Fatou component. Then $R_f$ is closed, reflexive and symmetric. Clearly, $\sim_f$ is the  minimal closed equivalence containing $R_f$. 

Let $R_f[x]$ be the fiber of $R_f$, which consists of all points $y$ related to $x$ under $R_f$. When $x\notin K$, $R_f[x]=\{x\}$. When $x\in K$, let $A(x)$ be the atom of $K$ containing $x$. Then $R_f[x]$ falls into three possibilities. First, if $x$ does not belong to the boundary of any bounded Fatou component, then $R_f[x]=A(x)$. Second, if $x$ lies in a bounded Fatou component  $W$, then $R_f[x]=A(x)\cup\overline{W}$. Finally, if there is a bounded Fatou component  $w$ with $x\in W$, then $R_f[x]=\overline{W}$. 

Clearly, $R_f[f(x)]=f\left(R_f[x]\right)$ holds for all $x\in\widehat{\mathbb{C}}$. Let $[z]_{\sim_f}$ be the class of $\sim_f$ containing $z$. Then $f\left([z]_{\sim_f}\right)=[f(z)]_{\sim_f}$ holds for all $z\in\widehat{\bbC}$, by the openness of $f:\widehat{\mathbb{C}}\rightarrow\widehat{\mathbb{C}}$ and the following.
\begin{lemm}[{\cite[Lemma 3.5]{LYY25}}]
Let  $g: X\rightarrow Y$ be an open map of a compactum $X$ onto another. Let $\mathcal{R}_X$ and $\mathcal{R}_Y$ be closed reflexive symmetric relations on $X$ and  on $Y$,  respectively, such that $g\left(\mathcal{R}_X[x]\right)=\mathcal{R}_Y[f(x)]$ for all $x\in X$. If $\sim_X$ (respectively, $\sim_Y$) is the minimal closed equivalence containing $\mathcal{R}_X$ (respectively, $\mathcal{R}_Y$), then every class $[x]_{\sim_X}$ of $\sim_X$ is sent by $g$ onto a class of $\sim_Y$. 	
\end{lemm}

This proves the first assertion. Below we continue to prove the second. 

To this end, we denote by  $\overline{f}:\Dc(f)\rightarrow\Dc(f)$  the induced map that sends $\delta$ to $f(\delta)$ and by $\pie: \widehat{\bbC}\rightarrow\Dc(f)$  the projection. Then $\pie$ is  monotone  and semi-conjugates $\left(\widehat{\bbC},f\right)$ with the induced system $\left(\Dc(f),\overline{f}\right)$. Clearly, $\Tc(f)=\pie(J)$, since by definition $\Tc(f)=\pie(K)$. Moreover, $\left(\Tc(f),\overline{f}\right)$ is a dendrite factor of $(K,f)$, provided that $K$ is connected.

Given $\delta\in\Dc(f)$ and $J_0$ the component of $J$ containing $\delta$, we will use the following separation theorem due to Janizewskiw to verify that $\delta\cap J_0$ is connected. This completes our proof.
\begin{lemm}[{\bf \cite[Theorem (5.1)]{UrsellYoung51}}]\label{lem:Janizewski}
If the common part of two continua $K,L\subset\widehat{\bbC}$ is disconnected  there are two points separated by $K\cup L$ but not separated by either $K$ or $L$.
\end{lemm}

In deed, if on the contrary $\delta\cap J_0$ were disconnected there would exist two points $x\ne y$ lying in $\widehat{\bbC}\setminus(\delta\cup J_0)$ that are separated by $\delta\cup J_0$, while neither 
$\delta$ nor $J_0$ separates $x$ from $y$. Let $U_x$ and $U_y$ be the components of $\widehat{\bbC}\setminus(\delta\cup J_0)$ that contain $x$ and $y$, respectively. Then one of them, say $U_x$, is contained in a bounded Fatou component, which is in turn contained in a single element of $\Dc(f)$. This is absurd, since the closure $\overline{U_x}$ intersects $\delta$ and $\delta'\cap\delta=\emptyset$ for all $\delta'\in\Dc(f)\setminus\{\delta\}$. 
\end{proof}


It follows that $\Dc^J(f)=\{\delta\cap J: \delta\in\Dc(f)\}$ is an USC decomposition of $J$ into sub-continua. Moreover, $f(\delta\cap J)=f(\delta)\cap J$ for all $\delta\in\Dc(f)$. Therefore, one can easily infer the following.
\begin{coro}\label{cor:DF_(J,f)}
If $J$ is connected then $\left(\Tc(f),\overline{f}\right)$ is a dendrite factor of $(J,f)$.
\end{coro}

\begin{rema}With Corollary \ref{cor:DF_(J,f)} we see that, to prove Theorem \ref{theo:MDF}, it suffices to show that for any dendrite factor $(T,g)$ the decomposition $\left\{\varphi^{-1}(u): u\in T\right\}$ given by the underlying semi-conjugacy $\varphi:J\rightarrow T$ is refined by $\Dc^J(f)$. We will do this in Theorem \ref{theo:C_for_core}. 
Also note that $\left(\Tc(f),\overline{f}\right)$ is generally not the maximal dendrite factor of $(K,f)$. For instance, if $f(z)=z^2$ then $K=\overline{\bbD}$ and $\Tc(f)$ is a singleton while $(K,f)$ has a dendrite factor which is topologically conjugate with $\left([0,1],x\mapsto x^2\right)$, via the semi-conjugacy $\varrho(z)=|z|$.
\end{rema}

Given a continuum $X\subset\widehat{\bbC}\setminus\{\infty\}$, let $\textbf{TH}(X)$  consist of $X$ and the components of
$\widehat{\bbC}\setminus X$ that do not contain $\infty$. Following \cite[Definition 14]{BCO11}, we call $\textbf{TH}(X)$ the {\bf topological hull} of  $X$. A compactum $X\subset\widehat{\bbC}$  is said to be {\bf full} if $\widehat{\bbC}\setminus X$ is connected or equivalently,  if $X=\textbf{TH}(X)$. Given a compactum $X\subset\widehat{\bbC}$,  we can identify an USC decomposition  $\Dc_X^\bbC$ whose elements are just the  singletons disjoint from $X$ and the topological hulls $\textbf{TH}(\delta)$ for atoms  $\delta$ of $X$.

By Moore's Decomposition Theorem \cite[Theorem 22]{Moore25}, we may pick a homeomorphism $h_K$ between the quotient space of $\Dc_K^\bbC$, still denoted by  $\Dc_K^\bbC$,  and $\widehat{\bbC}$ that sends $\{\infty\}\in\Dc_K^\bbC$ to $\infty\in\widehat{\bbC}$.  Let $m_K$ be the composite of  the projection $\widehat{\bbC}\rightarrow\Dc_K^\bbC$ followed by $h_K$. By \cite[Theorem 7]{LLY19}, the following is immediate.
\begin{lemm}\label{lem:m_K}
For any compactum $K$, $m_K$ is the finest monotone map of $\widehat{\bbC}$ onto itself that sends   $K$ to a Peano compactum. Moreover, if $K$ is full, so is $m_K(K)$.
\end{lemm}
When $K$ is the filled  Julia set of a polynomial, $m_K$ is just the monotone map $m_f$ obtained in \cite[Lemma 15]{BCO11} (for connected $K$) and in \cite[Theorem 19]{BCO13} (for disconnected $K$).  Notice that $m_K(K)$ may have interior points and one of its non-degenerate components is not a dendrite.

Hereafter in this subsection and in later subsections, we always assume that the underlying polynomial $f$  have a connected filled Julia set $K$. The B\"ottcher map $\phi_f$ for such $f$ is chosen to be a conformal mapping from $\widehat{\bbC}\setminus K$  onto
$\bbD^*=\left\{z\in\widehat{\bbC}: |z|>1\right\}$, fixing $\infty$,
such that $\phi_f\circ f(z)=\left(\phi_f(z)\right)^d$ holds for all $z\notin K$. For any $\theta\in\bbT=\bbR/\bbZ$, the ray connecting $\xi=e^{2\pi\textbf{i}\theta}$ to $\infty$ is mapped by  $\phi_f^{-1}$ to a smooth curve $R_K(\theta)$, to be called the {\bf external ray of $K$} (or the dynamic ray of $f$) at angle $\theta$.
Corresponding to $\xi$  there is a prime end of $\widehat{\bbC}\setminus K$ whose impression will be denoted by
${\rm Imp}(K,\xi)$ or by ${\rm Imp}(K,\theta)$, if we want to emphasize the angle $\theta$.

Before going on to prove Theorem \ref{theo:MDF}, we introduce further terminology below.
\begin{defi}\label{M_f^J}
Let  $\mathfrak{M}_f^J$ be the family of  all the USC decompositions $\Dc$ of $J$ into sub-continua such that (1) the resulting quotient space is a dendrite and (2)  $f(\delta)\in\Dc$ holds for all $\delta\in\Dc$.
\end{defi}

Clearly, the induced map  $\delta\xrightarrow{\widehat{f}}f(\delta)$ on the quotient space is well defined. Moreover, the natural projection $J\xrightarrow{\widehat{\pi}}\Dc$ is monotone and semi-conjugates $(J,f)$ with $\left(\Dc,\widehat{f}\right)$. In other words,  the induced system $\left(\Dc,\widehat{f}\right)$ is a dendrite factor of $(J,f)$. On the other hand, if $(T,g)$ is a dendrite factor of $(J,f)$ and $\varphi: J\rightarrow T$ the underlying semi-conjugacy, then all the point inverses $\varphi^{-1}(u), u\in T$ form an USC decomposition belonging to $\mathfrak{M}_f^J$.

Let $\Dc(f)$ be  given as in Definition \ref{def:sim_f} and $\Dc^J(f)=\left\{\delta\cap J: \delta\in\Dc(f)\right\}$.  Then $\Dc^J(f)$ belongs to $\mathfrak{M}_f^J$ by Theorem \ref{theo:invariance_sim}. Moreover,  Theorem \ref{theo:MDF} is resolved by the following.

\begin{main-theorem}\label{theo:C_for_core}
$\Dc^J(f)$ is the finest member of $\mathfrak{M}_f^J$.
\end{main-theorem}

\begin{proof} If $K=J$ then $\Dc^J(f)=\Dc_J^{PC}$ hence it refines $\Dc$. So we only consider the case $J\ne K$.

Given $\Dc\in \mathfrak{M}_f^J$. Let $\widehat{f}$ be the induced map and $\widehat{\pi}:J\rightarrow\Dc$ the natural projection. Since the resulting quotient space is a dendrite, every atom of $J$ is contained in a single element of $\Dc$.
Let us  fix a decomposition $\Dc^{\bbC}$ with two properties. First, every singleton contained in  $\widehat{\bbC}\setminus K$ is an element of $\Dc^{\bbC}$. Second, a bounded Fatou component $W$ is contained in an element of $\Dc$ if and only if $W$ is a bounded component of $\widehat{\bbC}\setminus \delta_W$ for some $\delta_W\in\Dc$. Moreover, when such  $\delta_W$ does not exist, every singleton lying in $W$ is an element of $\Dc^{\bbC}$.

Recall that $m_K\left(\overline{U_\infty}\right)$ is a non-degenerate Peano continuum with no cut point, where $U_\infty=\widehat{\bbC}\setminus K$. Still denote by $\widehat{\pi}:\widehat{\bbC}\rightarrow\Dc^\bbC$ the natural projection of $\Dc^\bbC$. Then there is a unique monotone map $\varphi:\widehat{\bbC}\rightarrow\widehat{\bbC}$ satisfying $\widehat{\pi}=\varphi\circ m_K$.
Moreover, it suffices to show that  $\varphi\circ m_K(W)$ is a singleton for all bounded Fatou components $W$ that are periodic.
Here, we only need to consider those  $W$ such that $m_K(W)$ is not a singleton. Fix such a component $W$. Then  $m_K(W)$ is a component of $\widehat{\bbC}\setminus m_K\left(\overline{U_\infty}\right)$ and, by \cite[p.126, Lemma 2]{Whyburn79}, the boundary of $m_K(W)$ is  a Jordan curve, to be denoted by $\Gamma_W$.

{\bf Claim 1}.  $m_K(\delta)\cap \Gamma_W$ is connected for each $\delta\in \Dc^\bbC$.

Suppose on the contrary there were $\delta\in\Dc^\bbC$ with $m_K(\delta)\cap \Gamma_W$ disconnected. Fix a separation $m_K(\delta)\cap \Gamma_W=E\cup F$. We can choose an open arc $\gamma\subset m_K(\delta)\setminus \Gamma_W$ connecting $x\in E$ to $y\in F$, since $m_K(\delta)$ is a Peano continuum.
As $\gamma\cup \Gamma_W$ is a $\theta$-curve, by \cite[p.123, $\theta$-curve Theorem]{Whyburn79} we know that $\widehat{\bbC}\setminus(\gamma\cup \Gamma_W)$ has two bounded components. One is $m_K(W)$. Let $U'$ be the other. Then $\partial U'$ intersects $\Gamma_W$ at a non-degenerate arc which, except for the two end points, is disjoint from the boundary   of $m_K(U_\infty)$. This is absurd, since $\partial W\subset\partial U_\infty$ and hence $\Gamma_W$ is contained in the boundary of $m_K\left(U_\infty\right)$.

{\bf Claim 2}.   $\varphi(\Gamma_W)$ is a singleton. Thus $\varphi\circ m_K(W)=\widehat{\pi}\left(\overline{W}\right)=\varphi(\Gamma_W)$ and we are done.

Suppose on the contrary that there were  $\delta'\ne\delta''$ in  $\Dc^\bbC$,  both intersecting $\partial W$, such that $\varphi\circ m_K(\delta')\ne\varphi\circ m_K(\delta'')$. Then $m_K(\delta')$ and $m_K(\delta'')$ would be disjoint. By {\bf Claim 1},  $A=m_K(\delta')\cap \Gamma_W$ and  $B=m_K(\delta'')\cap \Gamma_W$ are both connected. Since $\Gamma_W$ is a simple closed curve, we can find two open arcs $\alpha,\beta$ in $\Gamma_W\setminus(A\cup B)$  each of which connects a point on $A$ to a point on $B$. By {\bf Claim 1}, there is no $\delta\in \Dc^\bbC$ such that $m_K(\delta)$ intersects $\alpha$ and $\beta$ both. Thus $\varphi(\alpha)$ and $\varphi(\beta)$ are disjoint connected sets in $\widehat{\pi}(J)=\varphi\circ m_K(J)$. 
Here 
$\varphi(\alpha)\cup\left\{\varphi\circ m_K(\delta'),\varphi\circ m_K(\delta'')\right\}$ and $\varphi(\beta)\cup\cup\left\{\varphi\circ m_K(\delta'),\varphi\circ m_K(\delta'')\right\}$ are locally connected continua of 
$\widehat{\pi}(J)$ and they intersect at exactly the two points $\varphi\circ m_K(\delta')\ne \varphi\circ m_K(\delta'')$. This is absurd, since $\widehat{\pi}(J)$ is a dendrite.
\end{proof}

\subsection{Two Laminations and A Dichotomy of MDF's}\label{section:doubling_dynamics}
We shall address how $\left(\Tc(f),\overline{f}\right)$ is related to the dynamics of $\Big(\bbS^1,\sigma_d\Big)$. Some of those relations  are  useful  here as well as in continued projects. Let $\lambda_\bbR(f)$ be the closed equivalence on $\bbS^1$ that identifies $\xi_1,\xi_2$  if and only if the impressions ${\rm Imp}(K,\xi_i)$ both lie in the same atom of $K$. 

\begin{defi}\label{def:real_lamination}
We call $\lambda_\bbR(f)$ the {\bf real lamination} and $\lambda_\bbQ(f)=\lambda_\bbR(f)\cap\left\{(e^{2\pi{\mathbf i}\theta_1},e^{2\pi{\mathbf i}\theta_2}): \theta_i\in\bbQ\right\}$ the {\bf rational lamination}  of $f$.
\end{defi}

Kiwi's $\bbR$eal lamination is a closed equivalence on $\bbS^1$ satisfying the conditions (R1)-(R5) given in \cite{Kiwi04}. The above lamination $\lambda_\bbR(f)$ is more general, since we do not assume (R2), which require all the classes to be finite.  When $f$ has no irrationally neutral cycle, $\lambda_\bbR(f)$ satisfies (R1)-(R5) 
hence coincides with Kiwi's $\bbR$eal lamination. 
Notice that in such a case  we may use \cite[Theorems 2 and 3]{Kiwi04} to infer the equality $\lambda_\bbR(f)=\overline{\lambda_\bbQ(f)}$. 

\begin{defi}\label{def:dendrite_lamination}
Let $\tau(\xi)=\pie({\rm Imp}(K,\xi))$ for $\xi\in\bbS^1$. Let $\lambda_{\rm den}(f)$ be the equivalence on  $\bbS^1$ whose classes are  
$\left\{\tau^{-1}(u): u\in\Tc(f)\right\}$. We  call $\lambda_{\rm den}(f)$ the {\bf dendrite lamination} of $f$.
\end{defi}

The class of $\lambda_\bbR(f)$ and that of $\lambda_{\rm den}(f)$ containing $\xi$ will be denoted by $[\xi]_\bbR$ and $[\xi]_{\rm den}$, respectively. It is clear that $\lambda_\bbR(f)$ refines $\lambda_{\rm den}(f)$. If $K^o\ne\emptyset$ (thus $J\ne K$), $\lambda_{\rm den}(f)$ may have classes that are Cantor sets. In such a case, $\Tc(f)$ is different from the {\bf topological Julia set} given in \cite[Definition 1.6]{BBS22}. On the other hand, if $K=J$ then $\lambda_\bbR(f)=\lambda_{\rm den}(f)$. Let  the {\bf convex hull} of a set $X$, denoted by ${\rm CH}\left(X\right)$,  be the smallest closed convex set  containing $X$. The following is immediate.
\begin{lemm}\label{lem:usc_D}
 $\displaystyle \left\{{\rm CH}\left([\xi]_{\rm den}\right): \xi\in\bbS^1\right\}$ is a USC decomposition of the closed unit disk $\overline{\bbD}$.
\end{lemm}

Now we  show that $\Tc_1(f)$ (respectively, $\Tc_0(f)$) equals the minimal sub-continuum of $\Tc(f)$   containing $\pie(P_f)$ (respectively, $\pie\left({\rm Crit}(f)\cup P_f\right)$). 
To do that, we just obtain the following.

\begin{prop}
	\label{theo:CH_core_coro}
Let $Y\subset\Tc(f)$ be a compactum satisfying $\overline{f}(Y)\subset Y$ and $\pie(P_f)\subset Y$. Let $T_Y$ be the minimal  sub-continuum  of $\Tc(f)$  containing $Y$. Then  $\overline{f}\left(T_Y\right)\subset T_Y$.
\end{prop}

\begin{proof} Clearly, all the end points of $T_Y$ lie in $Y$. Since $\Tc(f)$ is a dendrite, for any  $y_1\ne y_2$ in $\Tc(f)$ there is an arc connecting $y_1$ to $y_2$, to be denoted by $[y_1,y_2]$. Since  $ \pie(P_f)\subset Y\subset T_Y$ by assumption, the inclusion $\overline{f}(u)\in T_Y$ holds for all $u\in\pie({\rm Crit}(f))$. Thus it suffices to show that $\overline{f}(u_0)\in T_Y$  for $u_0$ in $T_Y\setminus \pie({\rm Crit}(f))$. 

Let $Q_0$ be the component of $T_Y\setminus \pie({\rm Crit}(f))$ containing $u_0$. The closure $\overline{Q_0}$ is a dendrite, all of whose ends  belong to $Y$. Pick two ends $u_1, u_2$ of $\overline{Q_0}$ such that $u_0\in(u_1,u_2)$. Then $\overline{f}$ restricted to $[u_1,u_2]$ is injective. Thus $\overline{f}([u_1,u_2])$ is an arc and coincides with $\left[\overline{f}(u_1),\overline{f}(u_2)\right]$,  which is a subset of $T_Y$. This completes the proof, since $\overline{f}(u_0)\in \overline{f}([u_1,u_2])$.
\end{proof}

\begin{rema}\label{rem:lam_rel}
In our study, a lamination means a closed   equivalence on $\bbS^1$  that is {\bf tree-like}, so  that  any two classes $A\ne B$  are respectively contained in two disjoint arcs $I_A, I_B\subset\bbS^1$.  In other words, the convex hulls of two different classes under the hyperbolic metric  are always disjoint. Such a closed equivalence is called a {\em laminational equivalence relation} in \cite{ALRV20}, if in addition all the classes are required to be finite.    By definition,  both $\lambda_\bbR(f)$ and $\lambda_{\rm den}(f)$ are treelike, while they may have infinite classes. However, by \cite[Theorem 30]{BCO11} and Theorem \ref{theo:invariance_sim} we can infer that both $\lambda_\bbR(f)$ and $\lambda_{\rm den}(f)$ are  degree-$d$-invariant  equivalences in the sense of \cite[$\S2$]{MR4205644}.  Therefore, we will follow \cite{MR4205644} when introducing certain necessary notions, such as critical class and major.
\end{rema}

 \begin{defi}\label{def:critical_class}
We call $[\xi]_{\rm den}$  a {\bf critical class} of $\lambda_{\rm den}(f)$ provided that  $\sigma_d$ restricted to $[\xi]_{\rm den}$ is not injective. Similarly, we can define critical classes of $\lambda_\bbR(f)$.
\end{defi}

The above definition   is based on  \cite[Proposition 2.1]{MR4205644}. Routinely, we can verify that $[\xi]_{\rm den}$ is a critical class if and only if  there is at least one critical point of $f$ which is mapped by $\pie:\widehat{\bbC}\rightarrow\Dc(f)$ to $\tau(\xi)$.
Hereafter, if a critical class  consists of exactly two points, its convex hull (a geodesic arc) will be called a {\bf critical leaf}; if it contains at least three points, its convex hull will be called a {\bf critical gap}.

\begin{defi}\label{def:major}
The {\bf major} of $\lambda_{\rm den}(f)$ is the set of critical leaves and critical gaps.
\end{defi}

Let $\#A$ be the cardinality of a set. By using the semi-conjugacy $\tau:\bbS^1\rightarrow\Tc(f)$, we can infer the following dichotomy.

\begin{main-theorem}\label{theo:dichotomy}
If $f\in\Cc_d$ then either $h\left(\Tc(f),\overline{f}\right)=\log d$ or  $\#\Tc(f)=1$.
\end{main-theorem}

\begin{proof}
Let $\mu$ be the Haar measure of $\left(\bbS^1,\sigma_d\right)$. Then either $\mu\left(\tau^{-1}(u_0)\right)>0$ for some $u_0$ or $\mu\left(\tau^{-1}(u)\right)=0$ for all $u\in\Tc(f)$.

In the first case, $X=\tau^{-1}(u_0)$ is uncountable and there exists a smallest integer $n\ge0$ such that  $\sigma_d^{n}(X)=\sigma_d^{n+k}(X)$ for some $k\ge1$. Setting $Y=\sigma_d^n(X)\cup\cdots\cup\sigma_d^{n+k-1}(X)$, we can infer that $Y\subset\sigma_d^{-1}(Y)$ and that $\mu\left(Y\triangle\sigma_d^{-1}(Y)\right)=0$, where $Y\triangle\sigma_d^{-1}(Y)=\left[Y\setminus\sigma_d^{-1}(Y)\right]\cup\left[\sigma_d^{-1}(Y)\setminus Y\right].$
Since $\mu$ is ergodic with respect to $\sigma_d$, we also have $\mu\left(\bbS^1\setminus Y\right)=0$. It follows that $\tau^{-1}(u_0)=\bbS^1$, or equivalently, $\#\Tc(f)=1$.

In the second, from the tree-likeness of $\lambda_{\rm den}(f)$ one can infer that for all but countably many $u\in\Tc(f)$ the point inverse $\tau^{-1}(u)$ contains $\le2$ points. Thus  $\mu$ is supported on the union of $E_1(f)=\left\{w\in\bbS^1:\ \#\tau^{-1}(\tau(w))=1\right\}$ and $E_2(f)=\left\{w\in\bbS^1:\ \#\tau^{-1}(\tau(w))=2\right\}$. Set $\nu=\mu\circ\tau^{-1}$.  By \cite[formula (1.2)]{LedrappierWalters77}, we have
\begin{equation}
h_\mu(\sigma_d) 
=h_\nu\left(\overline{f}\right)+
\int_{\Tc(f)}h\left(\sigma_d,\tau^{-1}(y)\right)d\nu(y).
\end{equation}
Here  
$\displaystyle h\left(\sigma_d,X\right)=\lim\limits_{\varepsilon\rightarrow0}\limsup\limits_{n\rightarrow\infty}\frac{\log s_n(\varepsilon,X)}{n}$ 
for any compact $X\subset\bbS^1$, in which $s_n(\varepsilon,X)$ denotes the largest integer such that there exists an $(n,\varepsilon)$-separated subset $E\subset X$ with respect to $\sigma_2$. See \cite[p.169, Definitiosn 7.12 and 7.13]{Walters82}. Since $h\left(\sigma_d,\tau^{-1}(y)\right)=0$ for all but countably many $y\in\Tc(f)$, it is immediate that $h_\nu\left(\overline{f}\right)=h_\mu(\sigma_d)=\log d$. 
\end{proof}

\subsection{The Biaccessibility Dimension of Polynomials}\label{section:biac_dim} We will  define the biaccessibility dimension of $f$ and give a sufficient condition under which this dimension equals the ratio $\frac{h_{\rm GCE}(f)}{\log d}$.  

Hereafter, set $A(f)=\tau^{-1}\left(\Tc_0(f)\right)$ and denote by $\dim_HX$ the Hausdorff dimension  of $X$. Then $\dim_HA(f)\cdot\log d=h(A(f),\sigma_d)\ge h_{\rm GCE}(f)$.  This inequality becomes an equality under mild conditions. In deed, we have the following.
\begin{lemm}\label{lem:trivial_equality}
If  every cycle of $f$ is repelling, then $h_{\rm GCE}(f)=\dim_H A(f)\cdot\log d$.
\end{lemm}
\begin{proof}
Since we assume that every cycle of $f$ is repelling, the filled Julia set $K$ has no interior hence coincides with the Julia set $J$. By \cite[Theorem]{Kiwi04}, every class of $\lambda_\bbR(f)$ is a finite set. Moreover, we have $\lambda_\bbR(f)=\lambda_{\rm den}(f)$. This implies that $\tau:\bbS^1\rightarrow\Tc(f)$ is a finite-to-one map. By  \cite[Theorem 17]{Bowen71}, it is immediate that $h_{\rm GCE}(f)=h(A(f),\sigma_d)=\frac{\dim_H A(f)}{\log d}$. We are done.
\end{proof}

If $f$ has non-repelling cycles, $\dim_HA(f)\cdot\log d$ might be strictly  greater than $h_{\rm GCE}(f)$.
In such cases, we  consider a special subset $B(f)=\left\{\xi\in A(f): \#\tau^{-1}(\tau(\xi))=2\right\}$. Elements of $B(f)$ are called  {\bf biaccessibility angles} of $f$. 
\begin{defi}\label{def:B_f}
We call  $\dim_HB(f)$  the {\bf biaccessibility dimension} of $f$.
\end{defi}

The biaccessibility dimension $\dim_HB(f)$ is different from a similar notion (of the same name) from the literature that is by definition the Hausdorff dimension of a set containing $B(f)$. See for instance \cite[Appendix]{Dudko-Schleicher20} for such a notion limited to quadratics. The quantity $\dim_HB(f)$ is convenient in two ways. First, it applies to all polynomials. Second, under mild assumptions  we can relate it to $h_{\rm GCE}(f)$ as follows.
\begin{main-theorem}\label{theo:biac_dim}
If $\dim_H\Big(\overline{B(f)}\setminus B(f)\Big)=0$ then
$\displaystyle h_{\rm GCE}(f)=h\left(\overline{B\left(f\right)},\sigma_d\right)=\log d\cdot\dim_HB\left(f\right)$.\end{main-theorem}
\begin{proof}
The points $u$ in $\tau(A(f))\setminus \tau(B(f))$ are of two kinds. In the first category, the point inverse $\tau^{-1}(u)$ is a singleton. In the second, $\tau^{-1}(u)$ contains three or more points.   The points in the first category may be approximated by points in $\tau(B(f))$. The points in the second form a countable set. Thus $\tau\left(\overline{B(f)}\right)=\Tc_0(f)$ and $\left(\Tc_0(f),\overline{f}\right)$ is a factor of $\left(\overline{B(f)},\sigma_2\right)$.  By assumption, $\dim_H\Big(\overline{B(f)}\setminus B(f)\Big)=0$. So, every ergodic probability measure $\mu$ of $\left(\overline{B(f)},\sigma_2\right)$ with positive entropy is supported on $B(f)$. As $\tau$ restricted to $B(f)$ is two-to-one,  by  \cite[formula (1.2)]{LedrappierWalters77} we have $h_\mu(\sigma_2)=h_{\mu\circ\tau^{-1}}\left(\overline{f}\right)$. As $\mu$ is flexible, $\displaystyle h_{\rm GCE}(f)=h\left(\overline{B\left(f\right)},\sigma_2\right)$ by the variational principle. See for instance \cite[p.188, Theorem 8.6]{Walters82}. By \cite[Proposition III.1]{Furstenberg67}, this completes the proof.
\end{proof}

In the rest of this subsection, we give more lemmas  concerning  $\tau:\bbS^1\rightarrow\Tc(f)$
and $\left(\Tc_0(f),\overline{f}\right)$.  Given $T\subset \Tc(f)$, denote by $X_T$ the set of  all the convex hulls of  $\tau^{-1}(u)$ for $u\in T$. Moreover, a compact set $X\subset\overline{\bbD}$ is said to be {\bf saturated} (with respect to $\lambda_{\rm den}(f)$) provided that it contains the convex hull of every class $[\xi]_{\rm den}$ 
with $\xi$ running through $X\cap\bbS^1$.

\begin{lemm}\label{lem:convex_hull_converse}
If a compact set $X\subset\overline{\bbD}$ is saturated and convex, $\tau\left(X\cap\bbS^1\right)$ is a dendrite.
\end{lemm}
\begin{proof} Since $\tau\left(X\cap\bbS^1\right)$ is a compact subset of the dendrite $\Tc(f)$, it suffices to prove that $\tau\left(X\cap\bbS^1\right)$ is connected. We will use a proof by contradiction.

Suppose on the contrary that $\tau\left(X\cap\bbS^1\right)$ were not connected. Then there would be  $u_1\ne u_2$ in $\tau\left(X\cap\bbS^1\right)$ separated by 
a point $u$ in $\Tc(f)\setminus\{u_1,u_2\}$. Pick a separation $\Tc(f)\setminus\{u\}=T_1\cup T_2$, with $u_i\in T_i$. Then $\tau^{-1}(u)$ is a class of $\lambda_{\rm den}(f)$ whose convex hull separates $\tau^{-1}(u_1)$ from $\tau^{-1}(u_2)$ in $\overline{\bbD}$. This contradicts the convexity of $X$.
\end{proof}

\begin{lemm}\label{lem:convex_hull}
If $T\subset\Tc(f)$ is a dendrite,  $X_T$ is saturated and convex.
\end{lemm}
\begin{proof}
For any $u\in T$ let $X_u$ denote the convex hull of $\tau^{-1}(u)$. By definition,  $X_T=\bigcup_uX_u$. It will suffice to show that $X_T$ is closed and convex. To show that $X_T$ is closed, we may pick an arbitrary sequence $x_n\in X_T$ with $x_n\rightarrow x_\infty$ and find classes $[\xi]_{\rm den}$ and $[\xi_n]_{\rm den}\ (n\ge1)$ satisfying  $x_\infty\in{\rm CH}\left([\xi]_{\rm den}\right)$ and $x_n\in{\rm CH}\left([\xi_n]_{\rm den}\right)$. Set $u_n=\tau(x_n)$ and $u_\infty=\tau(x_\infty)$. Then $\{u_n\}\subset T$ and $u_n\rightarrow u_\infty$. Since $T$ is a dendrite, it follows that $u_\infty\in T$ and $x_\infty\in\tau^{-1}(u_\infty)\subset X$. 

Now we continue to prove the convexity of $X_T$ by contradiction. 
Suppose on the contrary that $X_T$ were NOT convex. Then we could find   $x_1\ne x_2$ in $X_T$ and a point $x$ in the interior of the line segment
$\overline{x_1x_2}$ from $x_1$ to $x_2$, such that $x\notin X_T$. As $X_T\cap\overline{x_1x_2}$ is compact and contains 
$\{x_1,x_2\}$, $\overline{x_1x_2}\setminus X_T$ has a unique component $\gamma$ (which is an arc) containing $x$. The two ends of $\gamma$ are contained respectively in the convex hull of $\tau^{-1}(t_1)$ and that of $\tau^{-1}(t_2)$ for some $t_1\ne t_2$ belonging to $T$. Respectively denote these convex hulls by $Y_1$ and $Y_2$. Pick  $t\in \Tc(f)$ that separates $t_1$ from $t_2$ in $\Tc(f)$. Then $t\in T$, since $T$ is assumed to be a dendrite. Let $Y$ denote the convex hull of $\tau^{-1}(t)$. Then $Y$ separates $Y_1$ from $Y_2$ in $\overline{\bbD}$, indicating that $Y\cap\gamma\ne\emptyset$. This is absurd, since $Y\subset X_T$ and $\gamma\cap X_T=\emptyset$.
\end{proof}

By combining Lemmas \ref{lem:convex_hull_converse} and \ref{lem:convex_hull}, we have the following.
\begin{prop}\label{prop:tau_and_dendrites}
For a compact set $X\subset\overline{\bbD}$ to be both convex and saturated, it is sufficient and necessary that $\tau\left(X\cap\bbS^1\right)$	be a dendrite.
\end{prop}

\subsection{On Real Polynomials and PCF Polynomials}\label{section:good_extension}
This subsection proves Theorem \ref{theo:good_extension}. To do that, we just need to obtain the following two theorems.
\begin{main-theorem}\label{theo:real}
For  $f\in\Cc_d$ whoise  critical points are all real,  the following both hold:
\begin{itemize}
\item[(1)] There is a closed interval $I\subset K\cap\bbR$ with $f(I)\subset I$  such that the restriction $\rho=\left.\pie\right|_I$ is monotone and semi-conjugates  $(I,f)$ to $\left(\Tc_0(f),\overline{f}\right)$.
\item[(2)] If $\nu$ is an ergodic probability measure of $(I,f)$ and $\nu(\{x\})=0$ for all $x\in I$, then  $\left(I,f,\nu\right)$ is isomorphic with $\left(\Tc_0(f),\overline{f},\nu\circ \varrho^{-1}\right)$ in the sense of \cite[p.58, Definition 2.4]{Walters82}.
\end{itemize}
\end{main-theorem}

\begin{main-theorem}\label{theo:PCF}
The statements below hold for each PCF polynomial $f\in\Cc_d$:
\begin{itemize}
\item[(1)]  The restriction $\varrho=\left.\pie\right|_{\Hc(f)}$ is monotone and semi-conjugates  $(\Hc(f),f)$ to $\left(\Tc_0(f),\overline{f}\right)$.
\item[(2)] If $\nu$ is an ergodic probability measure of $(\Hc(f),f)$ and $\nu(\{x\})=0$ for all $x\in\Hc(f)$, then  $\left(\Hc(f),f,\nu\right)$ is isomorphic with $\left(\Tc_0(f),\overline{f},\nu\circ \varrho^{-1}\right)$.
\end{itemize}
\end{main-theorem}

\begin{proof}[{\bf Proof for Theorem \ref{theo:real}}] Since all coefficients of $f$ are real,  $K$ is a continuum symmetric about the real axis $\bbR$. Using \cite[p.73, Boundary Bumping Theorem I]{Nadler92}, we can verify  that $K\cap\bbR$ is connected hence is an interval. Thus there exist $a<b$ with $[a,b]=K\cap\bbR$. Notice that $b$ and $a$ are just the landing points of the  rays at angle $0$ and  at angle $\frac12$.
Let $I\subset[a,b]$ be the minimal sub-interval that contains all the critical points of $f$ (since all of them are real). Then we immediately have $f(I)\subset I$, $\pie(I)=\Tc_0(f)$ and $h(I,f)=h(K\cap\bbR,f)$.

To establish item (1), we just need to show that $\rho=\left.\pie\right|_I$ is monotone.
To this end, we recall that each $\delta\in\Dc(f)$ is symmetric about $\bbR$. By \cite[p.73, Boundary Bumping Theorem I]{Nadler92}, we can infer that every $\delta\in\Dc(f)$ either is disjoint from $\bbR$ or intersects $\bbR$ at an interval.
Now it follows that $\left.\pie\right|_I$ is a monotone map of $I$ onto $\Tc_0(f)$, which  semi-conjugates $(I,f)$ with
$\left(\Tc_0(f),\overline{f}\right)$.

In the sequel, we deal with item (2).

As $K$ is symmetric about $\bbR$,  by \cite[p.73, Boundary Bumping Theorem I]{Nadler92}we may conclude that every bounded Fatou component $W$ either is disjoint from $\bbR$ or intersects $\bbR$ at an interval. If a Fatou component $W$ does intersect $\bbR$, then the restriction of $f$ to $\bigcup_{n\ge0}f^n\left(I\cap\overline{W}\right)$ has at most finitely many ergodic probability measures. Those measures are all supported on a cycle hence none of them equals $\nu$, since by assumption $\mu(\{x\})=0$ for all $x\in I$. Therefore, in order to establish item (2) we just  show that $\nu$ is supported on $X=\{z\in I: \rho^{-1}(\rho(z))\ \text{is\ a\ singleton}\}$. Indeed, for any $z\in I\setminus X$ let $\delta(z)$ be the unique element of $\Dc(f)$ containing $z$. Then $\delta(z)\cap I$ is a non-degenerate interval and there are two possibilities. First, $\delta(z)\cap I$ is wandering hence $\nu(\delta(z)\cap I)=0$. Second, it is eventually periodic hence   $\nu(\delta(z)\cap I)=0$ again.
We are done.
 \end{proof}


If $f$ is PCF and $K^o=\emptyset$, then $K=J$ and $\left.\pie\right|_{\Hc(f)}$ is a homeomorphism that conjugates $(\Hc(f),f)$ with $\left(\Tc_0(f),\overline{f}\right)$.
Therefore, to establish Theorem \ref{theo:PCF} we only consider the case   $K^o\ne\emptyset$.

Before going on to prove Theorem \ref{theo:PCF}, let us recall some notions and known facts about $\Hc(f)$. Following \cite{Poirier10},  by a {\bf regulated arc} we mean  a simple arc  $\gamma\subset K$ such that for every bounded Fatou component $W$ the common part $\gamma\cap\overline{W}$ is empty or a single point or a tree consisting of radial segments. Here a radial segment is defined as follows.

Every bounded Fatou component $W$ contains a unique point $c_W$ such that $f^k(c_W)$ is a critical point for some $k\ge0$. This point $c_W$ is referred to as the {\bf center} of $W$. Fix a conformal homeomorphism $\phi_W$ of $\mathbb{D}$ onto $W$ that sends the origin to $c_W$. A {\bf radial segment} then means the image under $\phi_W$ of a segment of the form $\{re^{2\pi{\mathbf i}\theta}: 0\le r\le 1\}$ for some $\theta\in\bbT$.

Since $\Tc_0(f)\subset \Tc(f)$ is the smallest subcontinuum containing $\pie({\rm Crit}(f))$, we immediately  have $\pie(\Hc(f))\supset\Tc_0(f)$.
Since $\pie^{-1}(\Tc_0(f))$ contains all vertices of $\Hc(f)$, we can obtain the reverse containment 
by showing that $\pie^{-1}(\Tc_0(f))\cap\Hc(f)$ is connected.
This is done in the following.

\begin{lemm}\label{2.1a}
$\rho=\left.\pie\right|_{\Hc(f)}$ is a monotone map of $\Hc(f)$ onto $\Tc_0(f)$.
\end{lemm}
\begin{proof}
We just need to show that $\pie^{-1}(M)\cap\Hc(f)$ is connected for any continuum $M\subset \Tc(f)$.
Suppose on the contrary that  $\pie^{-1}(M)\cap \Hc(f)$ were disconnected for a continuum $M\subset \Tc(f)$. By Lemma \ref{lem:Janizewski}, we can find two points $x\ne y$ that are separated by $\pie^{-1}(M)\cup \Hc(f)$ but neither of 
$\pie^{-1}(M)$ and $\Hc(f)$. Let $U_x$ and $U_y$ be the complementary components of  $\pie^{-1}(M)\cap \Hc(f)$ that contain $x$ and $y$, respectively. Then one of them, say $U_x$, is bounded and the boundary $\partial U_x$ intersects at least two components of $\pie^{-1}(M)\cap\Hc(f)$. It follows that $U_x$   is  contained in a bounded Fatou component $W$. This is absurd, since  $\overline{W}\cap\Hc(f)$ is a tree containing $\overline{U_x}\cap\Hc(f)$. \end{proof}

Now we are well prepared to prove Theorem \ref{theo:PCF}.

\begin{proof}[{\bf Proof for Theorem \ref{theo:PCF}}]
As $\pie\circ f(z)=\overline{f}\circ\pie(z)$ for all $z\in K$, Lemma \ref{2.1a} fixes  item (1). To deal with  item (2), we fix an arbitrary ergodic probability measure  $\nu$ of $(\Hc(f),f)$ such that  $\nu(\{x\})=0$ for all $x\in\Hc(f)$. Then, for any bounded Fatou component $W$, the restriction of $f$ to $\Hc(f)\cap\left(\bigcup_{n\ge0}f^n\left(\overline{W}\right)\right)$ has finitely many ergodic probability measures, each of which is supported on a cycle. Thus $\nu\left(\overline{W}\right)=0$ and  $\nu$ is supported on $X=\{z\in\Hc(f): \rho^{-1}(\rho(z))\ \text{is\ a\ singleton}\}$, since for any $z\in\Hc(f)\setminus X$ there is a bounded Fatou component $W$ with $\rho(z)=\pie(W)$. Here $X$ actually consists of all $z\in\Hc(f)$ with $\{z\}\in\Dc(f)$. Let $Y=\rho(X)$ and $\nu^*=\nu\circ \rho^{-1}$. Then $\rho\circ f(x)=\overline{f}\circ \rho(x)$ holds for all $x\in X_0$ and $\left.\rho\right|_{X}$ is invertible and conjugates $(X,f,\nu)$ with $(Y,\overline{f},\nu^*)$. In other words,  $f:\Hc(f)\rightarrow\Hc(f)$ (preserving $\nu$) and $\overline{f}:\Tc_0(f)\rightarrow\Tc_0(f)$ (preserving $\nu^*$) are isomorphic in the sense of \cite[p.58, Definition 2.4]{Walters82}. We are done.
\end{proof}


\subsection{On Renormalizable Polynomials}\label{section:small_J}
We consider renormalizable polynomials $f$ with a connected filled  Julia set $K$ such that  the restriction  $f_\#$ of $f^p(p\ge2)$  to an appropriate neighborhood $U$ of a continuum $Y\subset K$ is a polynomial-like map whose filled  Julia set is $Y$.

By the Straightening Theorem \cite[Theorem 1]{DH85}, there exist a unique polynomial $g$ up to affine conjugacy that is hybrid equivalent to 
$f_\#$, so that one can find  a Jordan domain $V$ with $\overline{V}\subset g(V)$  and a quasi-conformal mapping $\varphi: U\rightarrow V$ that conjugates $\left(U,f_\#\right)$ with $\left(V,g\right)$. 
\begin{defi}
We call $f_\#$ a {\bf renormalization} of $f$ and $g$ the {\bf straightening} of $f_\#$.
\end{defi}

Let $J_g$ be the Julia set of $g$ and $\pie_g:\widehat{\mathbb{C}}\rightarrow\Dc(g)$ the projection that semi-conjugates $\left(J_g,g\right)$ with its maximal dendrite factor $\left(\Tc\left(g\right),\overline{g}\right)$. Then $\left(\Tc\left(g\right),\overline{g}\right)$ is also the maximal dendrite factor of $\left(\partial Y,f_\#\right)$, via the monotone semi-conjugacy $\pie_g\circ\varphi$. Notice that   $\pie_g\circ\varphi$  semi-conjugates  $\left(Y,f_\#\right)$ with $\left(\Tc\left(g\right),\overline{g}\right)$, too.

Let $\pie: \widehat{\bbC}\rightarrow\Dc\left(f\right)$ be the natural projection and
$\left(\Tc(f),\overline{f}\right)$ the maximal dendrite factor of $(J,f)$, where $J=\partial K$. Since $J$ is also the Julia set of $f^p$ and since $\left(\Tc(f),\overline{f}^p\right)$ is the maximal dendrite factor of $(J,f^p)$, we have the following.
\begin{main-theorem}\label{small_J}
For any $\delta\in\Dc(f)$, either $\delta\cap Y$  is empty or $\pie_g\circ\varphi(\delta\cap Y)$ is a singleton. Consequently, the  dynamic core of $(J_g,g)$ is isomorphic with a subsystem of $\Big(\Tc_0(f),\overline{f}^p\Big)$.
\end{main-theorem}

We need two facts concerning how the atoms of $K$ are related to those of  $Y$.
First, given a component $C$ of $K\setminus Y$, by \cite[Corollary 4.10]{BOT17} and \cite[Lemma 3.5]{BOT17} we can infer that $\overline{C}\cap Y$ is contained in a single  impression of $\widehat{\bbC}\setminus Y$.
Second,  given  a compact set $X\subset\widehat{\bbC}$,  a component $U$ of $\widehat{\bbC}\setminus X$ and a compact set $L\subset\partial U$, by \cite[Theorem 1]{FLY22} we know that every atom of $L$ is contained in an atom of $X$.
From these results, it is easy to infer the following.

\begin{lemm}\label{lem:BOT17_key}
If $\delta$ is an atom of $K$ then $\delta\cap Y$ is either empty or an atom of $Y$.
\end{lemm}

\begin{proof}[{\bf Proof of Theorem \ref{small_J}}]
Recall that $\pie_g\circ\varphi$ restricted to $Y$ is a monotone map onto $\Tc\left(g\right)$ such that the boundary of any component $W$ of $Y^o$ (if $Y$ contains interior points) is mapped to a single point, denoted by $\delta_W$.
Applying Lemma \ref{lem:BOT17_key}, we can extend $\pie_g\circ\varphi: Y\rightarrow \Tc\left(g\right)$ to be a monotone map of $K$ onto $\Tc\left(g\right)$,  denoted by $m_Y$. Then, every component $C$ of $K\setminus Y$ is sent by $m_Y$ to the unique element $\delta_C\in \Tc\left(g\right)$ whose inverse under  $\pie_g\circ\varphi$ contains $\overline{C}\cap Y$. Concerning how to verify the continuity of $m_Y$, we refer to the proof of \cite[Theorem A]{BOT17}.

The point inverses of $m_Y$ now form an USC decomposition of $K$ such that no atom of $K$ and no interior component of $K$ is split. By Definition \ref{def:sim_f}, this decomposition is refined by $\Dc(f)$. Therefore, every  $\delta\in\Dc(f)$ either satisfies $\delta\cap Y=\emptyset$ or is such that $\delta\cap Y$ is a point inverse of $\pie_g\circ\varphi: Y\rightarrow \Tc\left(g\right)$. This completes our proof.\end{proof}

By Theorem \ref{small_J} we can infer that $\left(\Tc\left(g\right),\overline{g}\right)$ is conjugate with  a subsystem  of $\left(\Tc(f^p),\overline{f^p}\right)$.   Note that $\Tc(f^p)=\Tc(f)$ and $\overline{f^p}$ coincides with the $p$-th iterate $\overline{f}^p$ of $\overline{f}$. Moreover, the irreducible dynamic core $\left(\Tc_{\rm irr}(f^p),\overline{f^p}\right)$ of $\left(J,f^p\right)$ coincides with $\left(\Tc_{\rm irr}(f),\overline{f}^p\right)$, which is just the $n$-th iterated system of the irreducible dynamic core $\left(\Tc_{\rm irr}(f),\overline{f}\right)$ of $(J,f)$. Also note that every critical point of $f_\#$ is sent by $\varphi$ to a critical point of $g$. It follows that the irreducible dynamical core $\left(\Tc_{\rm irr}(g),\overline{g}\right)$ 
is conjugate with a subsystem of $\left(\Tc_{\rm irr}(f^p),\overline{f^p}\right)$. 

Since $\Tc_{\rm irr}(f)$ contains the measure center of $\left(\Tc_0(f),\overline{f}\right)$,   the following is immediate.
\begin{coro}\label{renormalizable_f}
$h_{\rm GCE}(f)\ge\frac1p h_{\rm GCE}(g)$.
\end{coro}

\section{Quadratic Theory of GCE}\label{sec:quadratic}
This section focuses on $f_c(z)=z^2+c$ with $c\in\M$, whose filled Julia set is denoted by  $K_c$  and  the Julia set by $J_c$. The materials of the subsections are outlined below.

Subsections \ref{section:major_H} and  \ref{section:atom} respectively prove the constancy of $h_{\rm GCE}(f_c)$  on two special types of sub-continua of $\M$. The first category includes the closure of the major cardioid $H_\heartsuit$ and hence the closure of any hyperbolic component, by Theorem \ref{theo:tuning_quadratic} and Lemma \ref{lem:monotone_tuning}. The second includes the prime end impressions of $\widehat{\bbC}\setminus\M$ and hence all the atoms of $\M$. 
These discussions are relevant due to two reasons. On the one hand, they help us to identify a continuous extension of the map $c\mapsto h(\Hc(f_c),f_c)$, defined for PCF quadratics $f_c$, over the whole Mandelbrot set. On the other, we can naturally induce a map of the locally connected model of $\M$ \cite{Douady93} which is comparable  with the map $\theta\mapsto\log r_\theta$ given in \cite[Theorem 8.2]{Tiozzo16}.  Recall that the MLC conjecture is open. So,  $\M$ might have non-degenerate atoms and even atoms having a bounded complementary component (which is then a queer component). 
Therefore, such an induced map (of the locally connected model of $\M$) is far from clear, although one may use the prime end theory to induce a map of the locally connected model of $\partial\M$, based on the previous map $\theta\mapsto\log r_\theta$.

Subsection \ref{sec:biac_dim_monotonicity} compares $h_{\rm GCE}\left(f_{c'}\right)$ to $h_{\rm GCE}\left(f_{c''}\right)$ for certain choices of $c'\ne c''$ in $\M$, aiming to illustrate some monotonicity of $h_{\rm GCE}:\M\rightarrow\bbR$.  

Subsection \ref{section:tuning} then shows that the dynamical core of $f_{c_H}$ is a dendrite factor of that of $f_{c_H\perp c}$, where $c_H\ne0$ is the center of a hyperbolic component $H$ and $c_H\perp c$ denotes the parameter obtained by tuning $c_H$ with $c\in\M\setminus\{\frac14\}$. With this, we may compare $h_{\rm GCE}\left(f_{c_H\perp c}\right)$ to $\max\left\{h_{\rm GCE}\left(f_{c_H}\right),\frac1p h_{\rm GCE}\left(f_{c}\right)\right\}$ and illustrate how the tuning operator introduced by Douady and Hubbard \cite{DH85b} plays a role in the GCE theory of quadratics.

Finally, subsection \ref{section:quadratic_theory} proves Theorem \ref{theo:properties_of_h}. 

\subsection{GCE is a Constant on the Closure of the Major Cardioid}\label{section:major_H}

Recall that the major cardioid  $H_\heartsuit$ of $\M$ consists of  all parameters $c$ such that $f_c(z)$ has an attracting fixed point and $\partial H_\heartsuit$ is a Jordan curve. We will obtain the following.

\begin{main-theorem}\label{theo:major_H}
If $c\in\overline{H_\heartsuit}$ then $\#\Tc(f_c)=1$ hence $h_{\rm GCE}\left(f_c\right)=0$.
\end{main-theorem}


\begin{proof}
For  $c\in H_{\heartsuit}\cup\{\frac14\}$, $J_c$ is a Jordan curve and $K_c$ is an element of the decomposition $\Dc(f_c)$. This indicates that $\#\Tc(f_c)=1$.   For $c\in\partial H_\heartsuit\setminus\{\frac14\}$, we separately deal with the three possible sub-cases that are mutually exclusive:  (1) $f_c(z)$ has a parabolic fixed point $z_0$, (2)  $f_c(z)$ has a Crem\'er fixed point $z_0$, (3)  $f_c(z)$ has a Siegel disk $U$.

In the first,  $f_c$  has a cycle of bounded Fatou components $V_j(1\le j\le k)$ with $z_0=\bigcap_i\partial V_i$. This implies that $\bigcup_{n=1}^\infty f_c^{-n}\left(\overline{V_1}\right)$ is a connected dense subset of $K_c$. By Definition \ref{def:sim_f}, we can infer that $\#\Tc(f_c)=1$.
In the second, $J_c$ has just one atom. See \cite[Example 3.6]{LYY25}. Thus $\Tc(f_c)$ is again a singleton.
In the last, the boundary  $\partial U$ is  either decomposable or indecomposable. 

 If $\partial U$ is decomposable, by \cite[Theorem 4.1]{GMO99} there exist a nowhere dense subcontinuum $B\subset J_c$ and a unique Cantor set $A\subset\bbT$ such that:
\begin{itemize}
\item[(a)] $B\supset\partial U$ and $f_c(B)=B$.
\item[(b)] Every prime end impression ${\rm Imp}(t)$ of $\widehat{\bbC}\setminus K_c$ with $t\in A$ intersects $\partial U$.
\item[(c)] Every $z\in B$ lies in some impression ${\rm Imp}(t)$ with $t\in A$.
\item[(d)] $B$ contains the critical orbit $\left\{f_c^n(0): n\ge0\right\}$.
\end{itemize}
Therefore, the continuum $B$ hence  $\left\{f_c^n(0): n\ge0\right\}$ is contained in a single element of  $\Dc(f_c)$, say $\delta(U)$,  which also contains $U$.
Recall that $f_c^{-1}\left(U\right)$ has exactly two components, one of which is $U$. Denote the other by $U'$. Since $\{0,c\}\subset\delta(U)$ and   $0\in\delta(U')$, we have $\delta(U')=\delta(U)=f_c^{-1}(\delta(U))$. This implies that $\delta(U)=K_c$ and that $\Tc(f_c)$ is a singleton.

If $\partial U$ is indecomposable, the previous sets $B\subset J$ and $A\subset\bbT$ are not known. And it is still unknown whether there is a quadratic $f_c(z)=z^2+c$ that has a fixed Siegel disk $U$ with $\partial U$ an indecomposable continuum. See Question 2 in section 8 of \cite{Rogers92-a}. This will be remedied by Proposition \ref{theo:confluence_lemma}, which then completes the proof.
\end{proof}

\begin{prop}\label{theo:confluence_lemma}
If $f_c(z)$ has an irrationally neutral fixed point $z_0$, then  $\#\Tc(f_c)=1$.
\end{prop}

We will need two results by Perez-Marco.

\begin{lemm}[{\rm\cite[Theorem 1]{Perez-Marco97}}]\label{Pe-1997-Thm1}
Let $\phi(z)=\exp(2\pi\textbf{i}\alpha)\cdot z+O(z^2)$, where $\alpha\in\bbR$, be a local holomorphic diffeomorphism. Let $U$ be a Jordan neighborhood of the indifferent fixed point $0$. Assume that $\phi$ and $\phi^{-1}$ are defined and univalent on a neighborhood of the closure of $U$. Then there exists a continuum  $X$ such that (i) $\widehat{\bbC}\setminus X$ is connected, (ii) $0\in X\subset\overline{U}$, (iii) $X\cap\partial U\ne\emptyset$, and (iv) $\phi(X)=X=\phi^{-1}(X)$. Moreover, if $\phi$ is not of finite order, $\phi$ is linearizable at $0$ if and only if $0\in\text{Int}(X)$.
\end{lemm}

\begin{lemm}[{\rm\cite[Theorem 2]{Perez-Marco97}}]\label{Pe-1997-Thm2}
Assuming the hypotheses in {\bf Lemma \ref{Pe-1997-Thm1}}, let $X$ be the Siegel compactum given by that theorem. Let $h:\widehat{\bbC}\setminus\overline{\bbD}\rightarrow\widehat{\bbC}\setminus X$ be a conformal representation of the exterior of $X$ such that $h(\infty)=\infty$. Then the map $g=h^{-1}\circ\phi\circ h$ extends to an analytic circle diffeomorphism of $\bbS^1=\partial\bbD$ with rotation number $\varrho(g)=\alpha$.
\end{lemm}

\begin{proof}[{\bf Proof for  Proposition \ref{theo:confluence_lemma}}] The projection $\pie: \widehat{\bbC}\rightarrow\Dc(f_c)$ has two special restrictions, each of which is a monotone map onto $\Tc(f_c)$. These monotone maps respectively semi-conjugate $(K_c,f_c)$ and $(J_c,f_c)$ with $\left(\Tc(f_c),\overline{f_c}\right)$, the maximal dendrite factor of $(J_c,f_c)$.

We only need to show that $\pie(z_0)=\pie(0)$, since in such a case $\pie(z_0)$ as a sub-continuum of $K_c$ is completely invariant and hence we have $K_c=\pie(z_0)$, or equivalently, $K_c\in\Dc(f_c)$. 

Suppose on the contrary that $\pie(z_0)\ne\pie(0)$. Then, they would be disjoint when considered as sub-continua of $K_c$. From this, we will infer a contradiction.

Assume that $f_c'(z_0)=\exp(2\pi\textbf{i}\alpha)$, where $\alpha\in\bbR\setminus\bbQ$.
Let $N$ be the union of $\pie(0)$ with all the external rays of $\widehat{\bbC}\setminus K_c$ whose limit sets intersect $\pie(0)$. Fix a Jordan domain $W_0$ which is contained in $\widehat{\bbC}\setminus N$ and contains $\pie(z_0)$. Then $f_c$ restricted to $W_0$ is univalent. Further fix $W\subset W_0$ with $\pie(z_0)\subset W$ and  $\overline{W}\subset(W_0\cap f(W_0))$. Then $f_c$ and $f_c^{-1}$ are both well defined and univalent on a neighborhood of the closure $\overline{W}$. In other words, $W$ and $\phi(z)=f_c(z+z_0)-z_0$ satisfy the hypotheses in Lemma \ref{Pe-1997-Thm1}. Pick a continuum $X$ that satisfies each of the following conditions:  (i) $\widehat{\bbC}\setminus X$ is connected, (ii) $0\in X\subset\overline{W}$, (iii) $X\cap\partial W\ne\emptyset$, and (iv) $\phi(X)=X=\phi^{-1}(X)$.

Further pick a conformal map $h:\widehat{\bbC}\setminus\overline{\bbD}\rightarrow\widehat{\bbC}\setminus X$ fixing $\infty$, as given in Lemma \ref{Pe-1997-Thm2}.  Then $g=h^{-1}\circ\phi\circ h$ is defined well on an open annulus $A$ such that  $\bbS^1$ is a component of $\partial A$. By Lemma \ref{Pe-1997-Thm2}, the map $g: A\rightarrow A$ extends to an analytic circle diffeomorphism of $\bbS^1$ with rotation number $\varrho(g)=\alpha$.
Note that $X+z_0\subset K$. This is clear since every $x\in X$ has a bounded orbit under $\phi$ hence every $z\in X+z_0$ has a bounded orbit under  $f_c$. We actually have $\partial X+z_0\subset J_c$. This containment  is clear, if $J_c=K_c$. If $J_c\varsubsetneq K_c$, or equivalently, if $f_c$ has an invariant Siegel disk $U$, one can infer the same containment in the following way.

Given a point $x\in U$ with $x-z_0\in X$. By the equality $X=\phi(X)$ we can infer that the orbit $\left\{f_c^n(x):n\ge0\right\}$ and hence its closure (a Jordan curve $J_x$) are contained in $X+z_0$. Let $W_x$ be the closed Jordan domain enclosed by $J_x$ that is contained in $K_c$. Then $W_x\subset U$, from which we can infer that the two sets $U\setminus W_x$ and $X+z_0$ have a common point $y$. Notice that the orbit $\left\{f_c^n(y):n\ge0\right\}$ and hence its closure (a Jordan curve $J_y$) are contained in $X+z_0$. Again, we can find a closed Jordan domain $W_y$, with $W_x\subset\text{Int}(W_y)\subset X+z_0$. Repeating this indefinitely, we can infer that $U\subset X+z_0$ and hence  $U\cap(\partial X+z_0)=\emptyset$. Since $X=\phi^{-1}(X)$, we can further infer that $\partial X$ is disjoint from all the backward iterates of $U$ that are contained in $X+z_0$.

By \cite[Theorem 1]{FLY22}, each atom of $\partial X+z_0$ is contained in a single atom of $J_c$. Therefore, each prime end impression of $\widehat{\bbC}\setminus(X+z_0)$ is  disjoint from or contained in $\pie(z_0)$. Fix an impression of $\widehat{\bbC}\setminus(X+z_0)$ contained in $\pie(z_0)$, which is just the cluster set $C(h,w)$ at some $w\in\bbS^1$. On the one hand, the equality $\phi(X)=X$ requires that the orbit of $C(h,w)$ under $\phi$ is entirely contained in $\pie(z_0)-z_0$ hence in $(\pie(z_0)-z_0)\cap\partial X$. On the other, the irrationality of the rotation number $\varrho(g)=\alpha$ requires that the orbit of $w$ under $g$ is dense in $\bbS^1$. Since $C(h,g(w))=\phi(C(h,w))$, we further infer that the orbit of $C(h,w)$ under $\phi$ is dense in $\partial X$. This is \textcolor{blue}{\bf absurd}, since we can use the fact $X\cap\partial U\ne\emptyset$ to infer that $(\pie(z_0)-z_0)\cap\partial X$ is a proper compact subset of $\partial X$.
\end{proof}

\begin{rema} 
By  Theorems \ref{theo:major_H} and \ref{theo:tuning_quadratic}, $h_{\rm GCE}(f_c)=0$ provided that $c$ can be connected to $0$ by finitely many hyperbolic components $H_1,\ldots,H_k$, so that $\bigcup_{i=1}^k\overline{H_i}$ contains both $c$ and $0$. By Corollary \ref{cor:GCE_level_0}, we even have $h_{\rm GCE}(f_c)=0$ for $c$ running through the central molecule $M_0$. 
\end{rema}

\subsection{GCE is a Constant on Each Atom of $\M$}\label{section:atom}
The local connectedness of $\M$ remains unknown, since the MLC conjecture is still open. Therefore, there might be a non-degenerate continuum $N$ that is an atom of $\partial\M$. This motivates us to address how $h_{\rm GCE}(f_c)$ varies when $c$  runs through  an atom of $\partial\M$. In deed, we  have the following.
\begin{main-theorem}\label{theo:atoms_of_M}
If $c',c''$ belong to the same atom of $\partial\M$ then $h_{\rm GCE}\left(f_{c'}\right)=h_{\rm GCE}\left(f_{c''}\right)$.
\end{main-theorem}

Let $\Dc_{\partial\M}^{PC}$ be the core decomposition of $\M$. Elements of $\Dc_{\partial\M}^{PC}$ are called atoms of $\partial\M$, each of which is also an atom of $\M$. See \cite[Lemmas 16 and 17]{BCO11} or  \cite[Theorem 7]{LLY19}. Note that every singleton in the interior of $\M$ is an atom of $\M$. Also note that there is a closed equivalence $\sim_\M$ on  $\bbS^1$ such that $w'\ne w''$ are equivalent if and only if the impressions ${\rm Imp}(\M,w')$ and ${\rm Imp}(\M,w'')$ are contained in the same atom of $\partial\M$.
If $w=e^{2\pi\textbf{i}\theta}$ we also write ${\rm Imp}(\M,\theta)$ instead of ${\rm Imp}(\M,w)$ and call it the prime end impression of $\widehat{\bbC}\setminus\M$ at angle $\theta\in\bbT$.
Combining this with \cite[Theorem 3]{Douady93}, 
we can infer that $c'\ne c''$ in $\partial\M$ are contained in the same atom of $\partial\M$  if and only if there exist two impressions ${\rm Imp}(\M,\theta')$ and ${\rm Imp}(\M,\theta'')$ whose union contains both $c'$ and $c''$. Therefore, with no loss of generality we may assume that $c', c''$ in Theorem \ref{theo:atoms_of_M} are contained in a single impression ${\rm Imp}(\M,\theta)$.

Since $f_c(z)=z^2+c$ has a neutral cycle if and only if $c$ lies on the boundary of a hyperbolic component $H\subset\M^o$, Theorem \ref{theo:atoms_of_M} shall be resolved by Lemmas \ref{lem:neutral_atom} and \ref{lem:constant_lam}.
\begin{lemm}\label{lem:neutral_atom}
If $c\in\partial H$ for some hyperbolic component $H$, then $\{c\}$ is an atom of $\partial\M$.
\end{lemm}

\begin{proof}
When $c=\frac14$ or when $H$ is not primitive, one can use known  results, such as those  recalled in \cite[Appendix D]{Tan00}, to conclude that $\{c\}$ is an atom of $\partial\M$. For  discussions in the language of fibers, see \cite{Schleicher04}. For a discussion in the language of atoms, see \cite{LY-2017}. So we only need to deal with  the case that $c\ne\frac14$ is the root of a primitive hyperbolic component.
In such a case, $f_c(z)$ has a cycle $z_1,\ldots,z_k$ of period $k>1$ such that the derivative of $f_c^k$ at every $z_i$ is equal to $1$.

Let $\Rc_c^{\pm}=\left\{\Phi^{-1}\left(re^{2\pi{\mathbf i}\theta^\pm}\right): r>1\right\}$ denote the two external rays of $\widehat{\mathbb{C}}\setminus\M$  that land at $c$. Then $\Gamma_c=\{c,\infty\}\cup\Rc_c^{+}\cup\Rc_c^{-}$ is a Jordan curve cutting $\M$ into two pieces. The one containing $H$ is called the wake of $\Gamma_c$, denoted by ${\rm wake}(c)$. The other contains $0$.
In the proof for \cite[Theorem 1.1(d)]{Tan00}, Tan concretely constructs a sequence of Jordan domains $P_n$ such that (1) every $P_n$ has a diameter $\le\frac{C}{n^2}$ for some constant $C>0$ and (2) $\partial P_n$ intersects $\M$ at exactly two points, one of which is $c$.
This implies that no point of $\M\setminus(\{c\}\cup{\rm wake}(c))$ lies in the atom of $\M$ that contains $c$, denoted by $\delta(c)$.
On the other hand, for any $c'\in\M\cap {\rm wake}(c)$ there is a Jordan curve $\Gamma$ disjoint from $\{c,c'\}$ intersecting $\partial\M$ at one or two points and separating $c$ from $c'$. This completes the proof, since it implies that no point $c'\in{\rm wake}(c)$ lies in the atom  $\delta(c)$.
\end{proof}

The Jordan curve $\Gamma_{c,c'}$ in the above proof  may be chosen in two ways, depending on whether $c'$ lies on $\overline{H}\setminus\{c\}$ or not. If $c'\in\overline{H}\setminus\{c\}$ then $\Gamma_{c,c'}$  consists of  two external rays of $\Phi$ landing at two points of $\partial H\setminus\{c,c'\}$, say $c_1$ and $c_2$, and the closure of an open arc inside $H$ connecting $c_1$ to $c_2$. If $c'$ does not lie on $\overline{H}\setminus\{c\}$ then it belongs to the wake of a parabolic parameter $c''\in\partial H$ (which is the root of a hyperbolic component) and we can choose $\Gamma_{c,c'}$ to be the union of $\{c'',\infty\}$ with the two external rays landing at $c''$.

\begin{lemm}\label{lem:constant_lam}
Let $c,c'\in{\rm Imp}(\M,\alpha)$, where $\alpha\in\bbT$ is chosen so that $\overline{H}\cap{\rm Imp}(\M,\alpha)=\emptyset$ for all hyperbolic components $H$. Then $\left(\Tc(f_c),\overline{f_c}\right)$ and  $\left(\Tc(f_{c'}),\overline{f_{c'}}\right)$ are topologically conjugate. Consequently,  $h_{\rm GCE}:\M\rightarrow\bbR$ is constant on ${\rm Imp}(\M,\alpha)$.
\end{lemm}

Before going on to prove Lemma \ref{lem:constant_lam}, we recall some terminology from \cite{Bandt-Keller92}. Given $\alpha\in\bbT$, let $X_\alpha=\left\{e^{2\pi\textbf{i}t}: t\ne\alpha\right\}$  and call  $\left\{\frac{\alpha}{2},\frac{1+\alpha}{2}\right\}$ a {\bf quadratic critical portrait}. See \cite[Definition 4.5]{Kiwi05}. For $t\in\bbT\setminus\{\alpha\}$,  let $l_j^\alpha\left(e^{2\pi\textbf{i}t}\right)=e^{\pi\textbf{i}(t+j+1)}\ (j=0,1)$. Then $l_0^\alpha,l_1^\alpha$ are maps of $X_\alpha$ into $\bbS^1$, which may be considered as the two ``inverse branches" of $\sigma_2(z)=z^2$. Here we note that $1\in l_1^\alpha(X_\alpha)$.

\begin{defi}\label{def:alpha_equiv}
By an {\bf $\alpha$-equivalence} we mean  a closed equivalence $\sim$ on $\bbS^1$ satisfying three conditions.: (a) 
$\displaystyle e^{\pi\textbf{i}\alpha}\sim e^{\pi\textbf{i}(1+\alpha)}$; (b)  $\displaystyle e^{2\pi\textbf{i}t}\sim e^{2\pi\textbf{i}s}$ implies $\displaystyle e^{2\pi\textbf{i}\cdot2t}\sim e^{2\pi\textbf{i}\cdot2s}$; (3) if $t,s\in\bbT\setminus\{\alpha\}$ then $e^{2\pi\textbf{i}t}\sim e^{2\pi\textbf{i}s}$ implies $l_0^\alpha\left(e^{2\pi\textbf{i}t}\right)\sim l_0^\alpha\left(e^{2\pi\textbf{i}s}\right)$ and $l_1^\alpha\left(e^{2\pi\textbf{i}t}\right)\sim l_1^\alpha\left(e^{2\pi\textbf{i}s}\right)$.
\end{defi}
The itinerary  of $\xi\in\bbS^1$ with respect to  $\left\{\frac{\alpha}{2},\frac{1+\alpha}{2}\right\}$, to be denoted by $I^\alpha(\xi)$, is defined to be an infinite sequence $(s_n)_{n=1}^\infty$ given by
\[
s_n=\left\{\begin{array}{ll}0&\text{if}\ \sigma_2^{n-1}(\xi)\in l_0^\alpha(X_\alpha),\\
1&\text{if}\ \sigma_2^{n-1}(\xi)\in l_1^\alpha(X_\alpha),\\
\star&\text{if}\  \sigma_2^{n-1}(\xi)\in\left\{e^{\pi\textbf{i}\alpha}, e^{\pi\textbf{i}(1+\alpha)}\right\}.\end{array}\right.
\]
In particular, the itinerary $I^\alpha(\alpha)$ is called the {\bf kneading sequence} of $\alpha$. Two $\alpha$-itineraries $I^\alpha(\xi')=(s_n)$ and $I^\alpha(\xi'')=(t_n)$ are said to be the same provided that  for all $n\ge1$ we have either $s_n=t_n$ or $s_n=\star$ or $t_n=\star$.
Note that the minimal closed equivalence  that identifies $\xi'\ne\xi''$ in $\bbS^1$ having the same $\alpha$-itinerary, to be denoted by $\approx_\alpha$, is an $\alpha$-equivalence. 
Also note that the intersection $\sim_\alpha$ of all the $\alpha$-equivalences as subsets of the product $\bbS^1\times\bbS^1$ is also an $\alpha$-equivalence. Clearly, $\sim_\alpha$ is a subset of $\approx_\alpha$.

\begin{defi}\label{def:minimal_alpha_equiv}
We call $\sim_\alpha$ the  {\bf minimal $\alpha$-equivalence} and $\approx_\alpha$ the {\bf dynamical $\alpha$-equivalence}. \end{defi}

\begin{rema}\label{rem:alpha_equiv_containments}
If $c\in\M$ is chosen so that the critical class of $\lambda_\bbR(f_c)$ contains $\frac{\alpha}{2}$, both $\lambda_\bbR(f_c)$ and $\lambda_{\rm den}(f_c)$ are $\alpha$-equivalences. Check Definitions \ref{def:real_lamination} and \ref{def:dendrite_lamination}. Routinely, one may verify the next two containments: (1) $\sim_\alpha \subset \lambda_\bbR\left(f_c\right)$; (2) $\approx_\alpha \subset \lambda_{\rm den}\left(f_c\right)$.
\end{rema}

\begin{proof}[{\bf Proof for Lemma \ref{lem:constant_lam}}] By assumption, $\overline{H}\cap{\rm Imp}(\M,\alpha)=\emptyset$ holds for each hyperbolic component $H$. So, $\alpha$ is irrational; moreover, $f_c(z)=z^2+c$ has no neutral cycle and hence $K_c=J_c$. Since $c'\in{\rm Imp}(\M,\alpha)$, $f_{c'}(z)=z^2+c'$ has no neutral cycle and $K_{c'}=J_{c'}$. Clearly, $\lambda_{\rm den}(f_{c})=\lambda_\bbR(f_c)$ and $\lambda_\bbR(f_{c'})=\lambda_{\rm den}(f_{c'})$.
By \cite[Theorem 1]{Kiwi04} and  Definition \ref{def:real_lamination},  the classes of $\lambda_\bbR(f_c)$ and those of $\lambda_\bbR(f_{c'})$ are finite. As $\alpha$ is irrational, we can use \cite[pp.193-194, Lemma 18.8]{Milnor06} to infer that the critical class of $\lambda_\bbR(f_c)$ (which contains $e^{\pi\textbf{i}\alpha}$) is not periodic. This implies that the kneading sequence $I^\alpha(\alpha)$ is not periodic and that $\sim_\alpha=\approx_\alpha$, by \cite[Theorem 1]{Bandt-Keller92}. By \cite[Theorem 3]{Kiwi05}, we further have $\lambda_\bbR(f_c)=\sim_\alpha=\lambda_\bbR(f_{c'})$. Therefore, the MDF  of $(J_c,f_c)$  and that of $(J_{c'},f_{c'})$ are topologically conjugate. So are $\left(\Tc_0(f_c),\overline{f_c}\right)$ and $\left(\Tc_0(f_{c'}),\overline{f_{c'}}\right)$.  This completes the proof.  \end{proof}



\subsection{The Monotonicity Theorem of GCE}\label{sec:biac_dim_monotonicity}
Let $A(f_c)=\tau^{-1}\left(\Tc_0(f_c)\right)$  and $B(f_c)$ the intersection of $A(f_c)$ with $\left\{\xi\in\bbS^1: \#\tau^{-1}(\tau(\xi))=2\right\}$. 
We will compare  $h_{\rm GCE}(f_c)$ with $h_{\rm GCE}(f_{c'})$ for $c\ne c'$ in $\M$ that satisfy certain conditions. See Theorem \ref{theo:monotonicity}.

Recall that every class $[w]_{\rm den}$ of $\lambda_{\rm den}(f_c)$ has a negative counterpart $[-w]_{\rm den}$ and that there is a unique critical class, which is just $\tau^{-1}\left(\pie(0)\right)$.  Moreover, the filled  Julia set $K_c$ is symmetric about the origin, so is the  major of $\lambda_{\rm den}(f_c)$. Let  $G_0(f_c)$ denote the convex hull of the critical class of $\lambda_{\rm den}(f_c)$ and 
 $G_1(f_c)$ the convex hull of  $\tau^{-1}\Big(\pie(c)\Big)$.
By Definition \ref{def:major}, $G_0(f_c)$  is the major of $\lambda_{\rm den}(f_c)$. Compare \cite[Definition 2.3]{BBS22}.  

If  $G_0(f_c)=G_1(f_c)=\overline{\bbD}$,  there is no minor at all and $\#\Tc(f_c)=1$. 
If $G_0(f_c)\ne G_1(f_c)$, $\overline{\bbD}\setminus\left(G_0(f_c)\cup G_1(f_c)\right)$ has a component  $Q(f_c)$ whose closure $\overline{Q(f_c)}$ intersects both $G_0(f_c)$ and $G_1(f_c)$.
We call $\overline{Q(f_c)}\cap G_1(f_c)$ the {\bf minor} of $\lambda_{\rm den}(f_c)$. 
If $\#\Tc(f_c)>1$, the minor of $\lambda_{\rm den}(f_c)$ is either a point  or a geodesic arc joining  $\xi'$ to $\xi''$. In the former case, let $\Ic(f_c)=\{\xi\}$; in the latter, let $\Ic(f_c)$ be  the subarc $[\xi',\xi'']$ of $\bbS^1\setminus Q(f_c)$.

\begin{defi}\label{minor_major} We call $\Ic(f_c)$ the {\bf characteristic arc} of $f_c$.  \end{defi}

When the minor of $\lambda_{\rm den}(f_c)$ is  a geodesic arc joining two points $\xi'\ne\xi''$, the difference $\bbS^1\setminus\sigma_2^{-1}\left(\{\xi'\,\xi\}\right)$ consists of four disjoint open arcs. Let $X_c$ be the union of the two that are intersect $G_0(f_c)$.  Then one can easily infer the following.

\begin{lemm}\label{lem:unlinking}
Given $\eta'\ne\eta''$ in $\bbS^1$, if for all $k\ge0$ the geodesic arc from $\sigma_2^k(\eta')$ to $\sigma_2^k(\eta'')$ is disjoint from 
 ${\rm CH}(X_c)$, then $[\eta']_{\rm den}$ contains $\eta''$ and  is disjoint from $\displaystyle\bigcup_{j\ge0}\sigma_2^{-j}\Big(G_0(f_c)\cap\bbS^1\Big)$.
\end{lemm}

Lemmas \ref{lem:non-critical_class} and \ref{lem:inclusion_lemma} are both useful.

\begin{lemm}\label{lem:non-critical_class}
Given $f_c(z)=z^2+c$ with $\#\Tc(f_c)>1$, if  a class $[\xi]_{\rm den}$ is  disjoint from the grand orbit of the critical class,  then $\#[\xi]_{\rm den}<\infty$.
\end{lemm}
\begin{proof}Let $\Lc(f_c)$ consist of three types of sets: (1) the singletons $\{\xi\}$ which are classes of $\lambda_{\rm den}(f_c)$, (2) the geodesic arcs  connecting two points which constitute a  class of $\lambda_{\rm den}(f_c)$,  and (3) the geodesic arcs connecting two points on the boundary of a gap which is the convex hull of a class of  $\lambda_{\rm den}(f_c)$. Then 
$\Lc(f_c)$ is a geometric lamination invariant  under $\sigma_2$, in the sense of \cite[Definition II.4.2]{Thurston09}.
We only consider classes 	$\#[\xi]_{\rm den}\ge3$ that contain at least three points. By \cite[II.5.2]{Thurston09} we know that $[\xi]_{\rm den}$ is eventually periodic under $\sigma_2$. Pick $k,p\ge1$ such that $\sigma_2^{k+p}\left([\xi]_{\rm den}\right)=\sigma_2^{k}\left([\xi]_{\rm den}\right)$, or equivalently, $\left[\sigma_2^{k+p}(\xi)\right]_{\rm den}=\Big[\sigma_2^{k}(\xi)\Big]_{\rm den}$. Since by assumption $[\xi]_{\rm den}$ is disjoint from the grand orbit of the critical class, $\sigma_2^p$ restricted to the class $\Big[\sigma_2^{k}(\xi)\Big]_{\rm den}$ is a homeomorphism onto itself. By \cite[pp.194-193, Lemma 18.8]{Milnor06}, 
we have $\#[\xi]_{\rm den}<\infty$.
\end{proof}

\begin{lemm}[{\bf Inclusion Lemma}]\label{lem:inclusion_lemma}
 Let  $c\ne c'$ in $\M$ be two parameters such that  $\Ic(f_{c})$ is a non-degenerate arc containing $\Ic(f_{c'})$. Then  $B(f_{c})\setminus B(f_{c'})$ is a countable set.
\end{lemm}
\begin{proof}  Since  $\Ic(f_{c})\supsetneq\Ic(f_{c'})$ by assumption,  the convex hull ${\rm CH}(X_c)$ given in Lemma \ref{lem:unlinking}  the critical gap $G_0(f_{c'})$ of $\lambda_{\rm den}(f_{c'})$.
Notice that, except for countably many choices of  $\xi\in B(f_c)$, the point $\tau(\xi)$ does not belong to the forward orbit of $\pie(0)$ under $\overline{f_c}$. Let $\xi\in B(f_c)$ be such a point and $\xi'\in B(f_c)$ the unique point with $\tau(\xi')=\tau(\xi)$. Then the forward orbit of $\xi$  under $\sigma_2(w)=w^2$ and that of $\xi'$ never intersect $\Ic(f_c)$. These two orbits are both disjoint from ${\rm CH}(X_c)$.
Let $\gamma_0$ be the minimal arc in $\Tc_0(f_c)$ that connects $\pie(0)$ to $\pie(c)$. Then $\Tc_0(f_c)=\bigcup_{j=0}^k\overline{f_c}^j(\gamma_0)$ for some $k\ge1$. It follows that  $\tau(\xi)=\overline{f_c}^j(u)$ for some interior point $u$  of $\gamma_0$ .  By the choice of $\xi$, we can infer that $\tau^{-1}(u)$ consists of exactly two points, say $\eta\ne \eta'$. Clearly, $\{\xi,\xi'\}=\sigma_2^j(\{\eta,\eta'\})$.

Recall that for all $j\ge0$ the geodesic arc $\Gamma_j$ connecting $\sigma_2^j(\eta)$ to $\sigma_2^j(\eta')$ is disjoint from $\overline{U_0}$. By Lemma \ref{lem:unlinking},  $\eta$ and $\eta'$ are contained in the same class  of $\lambda_{\rm den}(f_{c'})$, to be denoted by $[\eta]'_{\rm den}$, which  is disjoint from $\displaystyle\bigcup_{j\ge0}\sigma_2^{-j}\Big(G_0(f_c)\cap\bbS^1\Big)$. By  Lemma \ref{lem:non-critical_class},  we have $\#[\eta]'_{\rm den}<\infty$. It follows that either $\eta\in B(f_{c'})$ or $\#[\eta]'_{\rm den}$ is an integer no less than three.  Since $\lambda_{\rm den}(f_{c'})$ has at most countably many classes which have $\ge3$ but finitely many elements, this completes the proof.  \end{proof}

Now we are well prepared to obtain the following.

\begin{main-theorem}[{\bf Monotonicity Theorem}]\label{theo:monotonicity} 
If $\Ic(f_{c})\supset \Ic(f_{c'})$ then $h_{\rm GCE}(f_c)\le h_{\rm GCE}(f_{c'})$.
\end{main-theorem}
\begin{proof} 
With no loss of generality, we may assume that $\Ic(f_{c})\ne \Ic(f_{c'})$. Then  $\Ic(f_{c})$ is a non-degenerate arc and  $\Tc_0(f_c)$ is a tree. Moreover, the complement of any branch point $u_0$ in $\Tc_0(f_c)$ has finitely many components. This implies that $\overline{B(f_c)}\setminus B(f_c)$ is countable. By Theorem \ref{theo:biac_dim}, 
\begin{align}
h_{\rm GCE}(f_c)&=h\left(\overline{B(f_c)},\sigma_2\right)=\log2\cdot\dim_H\overline{B(f_c)}=\log2\cdot\dim_HB(f_c).
\end{align}
By Lemma \ref{lem:inclusion_lemma},  $B(f_{c})\setminus B(f_{c'})$ is also countable. There are two possible cases. First, every cycle of $f_{c'}$ is repelling hence we can apply Lemmas \ref{lem:trivial_equality} and \ref{lem:non-critical_class} to infer that 
\begin{align}
h_{\rm GCE}(f_c)&\le \log 2\cdot \dim_H\overline{B\left(f_{c'}\right)}\le \log 2\cdot \dim_HA\left(f_{c'}\right)=h_{\rm GCE}\left(f_{c'}\right).
\end{align}
Second, $f_{c'}$ has a non-repelling cycle hence $\Ic(f_{c})$ is a non-degenerate arc and  $\overline{B(f_c)}\setminus B(f_c)$ is  countable, too. Therefore, by Theorem \ref{theo:biac_dim} and the countability of $B(f_{c})\setminus B(f_{c'})$, we have
\begin{align}
h_{\rm GCE}(f_c)&=\log2\cdot\dim_HB(f_c)\le\log2\cdot\dim_HB(f_{c'})=h_{\rm GCE}(f_{c'}).
\end{align}
This completes the proof.\end{proof}

\begin{rema} 
If  $f_c, f_{c'}$ are PCF or parabolic and if  $c<c'$ (sometimes also written as $c\prec c'$) in the sense of \cite[$\S20.4$]{DH84}, then $\Ic(f_c)\supset\Ic(f_{c'})$. Therefore,  the first part of \cite[Lemma 2.5]{Dudko-Schleicher20} is extended to GCE by Theorem \ref{theo:monotonicity}. Note that  $\dim_HB(f_c)$ (see Definition \ref{def:B_f}) is different from $B_{\rm top}(c)$, since the latter is defined to be the Hausdorff dimension of a larger set containing $B(f_c)$. See \cite[Appendix]{Dudko-Schleicher20}. So, Lemma \ref{lem:inclusion_lemma} is a modified version of the second part of \cite[Lemma 2.5]{Dudko-Schleicher20}. 
 \end{rema}

\subsection{The Tuning Theorem of GCE}\label{section:tuning}
We will address how  $h_{\rm GCE}\left(f_{c_H\perp c}\right)$ is related to $h_{\rm GCE}\left(f_{c_H}\right)$ and $h_{\rm GCE}\left(f_{c}\right)$, where $H$ is a fixed hyperbolic component $H$ of period $p\ge2$ and $c_H$ its center. Notice that $c_H$  is a periodic point under $f_{c_H}(z)=z^2+c_H$. Due to \cite{DH85b}, it is known there is a continuous injection $\iota_H:\M\setminus\{\frac14\}\rightarrow\M$, satisfying $\iota_H(0)=c_H$, such that for any $c\in\M\setminus\frac14$  the $p$-th iterate  $f_{\iota_H(c)}^p$ has a quadratic-like restriction  which is hybrid equivalent to $f_c$. Moreover, if $H$ is primitive then $c$ may be chosen to be $\frac14$.
Set $c_H\perp c=\iota_H(c)$ and call  $\M_H=\Big\{c_H\perp c: c\in\M\setminus\{\frac14\}\Big\}$ the {\bf baby Mandelbrot set} of $H$. We have the following.

\begin{main-theorem}[{\bf Tuning Theorem}]\label{theo:tuning_quadratic}
$\left(\Tc_0\left(f_{c_H}\right),\overline{f_{c_H}}\right)$ is a dendrite factor of
$\left(\Tc_0\left(f_{c_H\perp c}\right),\overline{f_{c_H\perp c}}\right)$. Moreover,  $h_{\rm GCE}\left(f_{c_H\perp c}\right)=\max\left\{h_{\rm GCE}\left(f_{c_H}\right),\frac1p h_{\rm GCE}\left(f_{c}\right)\right\}$. 
\end{main-theorem}

Recall that ${\rm Imp}\left(K_c,\xi\right)$ denotes the  impression of $\widehat{\bbC}\setminus K_c$ at $\xi\in\bbS^1$
and that  $\tau^{c_H}:\bbS^1\rightarrow\Tc\left(f_{c_H}\right)$ (see Definition \ref{def:dendrite_lamination}) sends $\xi$ to  the element of $\Dc\left(f_{c_H}\right)$ that contains ${\rm Imp}\left(K_{c_H},\xi\right)$.  
By definition, $\tau^{c_H}$ and $\tau^{c_H\perp c}$ semi-conjugate $\left(\bbS^1,\sigma_2\right)$ with $\big(\Tc\left(f_{c_H}\right),\overline{f_{c_H}}\Big)$  and $\big(\Tc\left(f_{c_H\perp c}\right),\overline{f_{c_H\perp c}}\Big)$, respectively.

Let $W_0$ be the bounded Fatou component of $f_{c_H}$ containing $0$. The critical class of $\lambda_{\rm den}\left(f_{c_H}\right)$, denoted by $[w_0]_{\rm den}$, consists of all $\xi\in\bbS^1$ satisfying  ${\rm Imp}\left(K_{c_H},\xi\right)\subset\partial W_0$.
Clearly, we have $\tau^{c_H}\left([w_0]_{\rm den}\right)=\pie(0)=\pie(W_0)$.
As $f_{c_H}$  has no irrationally neutral cycle, by \cite[Theorems 1-3]{Kiwi04} we can infer that all the classes of $\lambda_\bbR\left(f_{c_H}\right)$ are finite sets and each of them is either contained in or disjoint from $\left\{e^{2\pi{\mathbf i}\theta}: \theta\in\bbQ\right\}$.
  Also note that each class of $\lambda_{\rm den}\left(f_{c_H}\right)$ that is not a class of $\lambda_\bbQ\left(f_{c_H}\right)$ allows two possibilities: (1) a finite  class of $\lambda_\bbR\left(f_{c_H}\right)$ disjoint from  $\left\{e^{2\pi{\mathbf i}\theta}: \theta\in\bbQ\right\}$  whose forward orbit never enters $[w_0]_{\rm den}$;
 (2) an uncountable set mapped by certain iterate of $\sigma_2$ onto $[w_0]_{\rm den}$.
 Similar results are known even in higher degrees \cite[Proposition 3.7(iii)-(iv)]{IK12}.
 
 Denote by $A\left(f_{c_H}\right)$  the pre-image of $\Tc_0\left(f_{c_H}\right)$ under $\tau^{c_H}$  and by $A\left(f_{c_H\perp c}\right)$
that of $\Tc_0(f_{c_H\perp c})$ under  $\tau^{c_H\perp c}$. Clearly, 
$h_{\rm GCE}\left(f_{c_H\perp c}\right)\le \dim_HA\left(f_{c_H\perp c}\right)\cdot\log2$.

\begin{proof}[{\bf Proof for Theorem \ref{theo:tuning_quadratic}}]
Let $\displaystyle X_H=\left\{\xi\in A\left(f_{c_H\perp c}\right): \sigma_2^k(\xi)\cap[w_0]_{\rm den}=\emptyset (k\ge0)\right\}$, where $[w_0]_{\rm den}$ is the critical class of $\lambda_{\rm den}\left(f_{c_H}\right)$. Then  every class 
 of $\lambda_{\rm den}\left(f_{c_H\perp c}\right)$ containing $\xi\in X_H$ is also  a class of $\lambda_{\rm den}\left(f_{c_H}\right)$. If  $\xi\in A\left(f_{c_H\perp c}\right)\setminus X_H$, then $f_{c_H}$ has a bounded Fatou component $W_\xi$, with ${\rm Imp}\left(K_{c_H},w\right)\subset\partial W_\xi$, such that the class of $\lambda_{\rm den}\left(f_{c_H\perp c}\right)$ containing $\xi$ is entirely contained in the class of $\lambda_{\rm den}\left(f_{c_H}\right)$ containing $\xi$, denoted by $[\xi]_{\rm den}$. Here we have $\tau^{c_H}\left([\xi]_{\rm den}\right)=\pie(W)$.  Moreover,
\begin{align}
\Tc\left(f_{c_H\perp c}\right)&=\left\{\tau^{c_H\perp c}(\xi): \xi\in A\left(f_{c_H\perp c}\right)\right\},
\\
\Tc\left(f_{c_H}\right)&=\Big\{\tau^{c_H}(\xi): \xi\in A\left(f_{c_H}\right)\Big\}.
\end{align}
Define $m: \Tc\left(f_{c_H\perp c}\right)\rightarrow \Tc\left(f_{c_H}\right)$ in the following way. If $\xi\in X_H$, set $m\left(\tau^{c_H\perp c}(\xi)\right)=\tau^{c_H}(\xi)$; otherwise, set $m\left(\tau^{c_H\perp c}(\xi)\right)=\pie(W_\xi)$. Notice that, in the latter case, we have $\pie(W_\xi)=\tau^{c_H}(\xi)$. 

Let $A_Y$ denote the preimage of $Y\subset\Tc\left(f_{c_H}\right)$ under $\tau^{c_H}$. Then $\tau^{c_H\perp c}$ is a closed mapping and $m^{-1}(Y)=\tau^{c_H\perp c}(A_Y)$ holds for any closed set $Y\subset\Tc\left(f_{c_H}\right)$. Such $Y$ is necessarily compact, so is $m^{-1}(Y)$. This resolves the continuity of 
$m: \Tc\left(f_{c_H\perp c}\right)\rightarrow \Tc\left(f_{c_H}\right)$. 
The surjection of $m$ comes from Lemma \ref{lem:convex_hull_converse}. Moreover, for any $w\in A\left(f_{c_H\perp c}\right)$  we have
\begin{equation}
m\circ \overline{f_{c_H\perp c}}\left(\tau^{c_H\perp c}(w)\right)=m\left(\tau^{c_H\perp c}(\sigma_2(w))\right)=
\tau^{c_H}(\sigma_2(w))=\overline{f_{c_H}}\circ\tau^{c_H}(w).
\end{equation}
Thus $m$ semi-conjugates $\left(\Tc_0\left(f_{c_H\perp c}\right),\overline{f_{c_H\perp c}}\right)$ with $\left(\Tc_0\left(f_{c_H}\right),\overline{f_{c_H}}\right)$. 
Let $\rho$ be the restriction of $\pie: \widehat{\bbC}\rightarrow\Dc\left(f_{c_H}\right)$ to $\Hc\left(f_{c_H}\right)$.  
\begin{figure}[ht]
\vspace{-0.05cm}
\begin{tabular}{ccc}
\begin{tikzpicture}
\matrix (m) [matrix of math nodes,row sep=2.75em, column sep=6em,minimum width=2em]
  {
   \Hc\left(f_{c_H}\right)  &  \Hc\left(f_{c_H}\right)  \\   \Tc_0\left(f_{c_H}\right) & \Tc_0\left(f_{c_H}\right)\\};
  \path[-stealth]
    (m-1-1) edge [double] node [right] {$\rho$}
    (m-2-1) edge  node [below] {$f_{c_H}$} (m-1-2) 
    (m-2-1.east|-m-2-2) edge node [above] {$\overline{f_{c_H}}$} (m-2-2)
    (m-1-2) edge [double] node [left] {$\rho$}  (m-2-2);
\end{tikzpicture}
&\hspace{1cm}&
\begin{tikzpicture}
\matrix (m) [matrix of math nodes,row sep=2.75em, column sep=6em,minimum width=2em]
  {
   \Tc_0\left(f_{c_H\perp c}\right)   & \Tc_0\left(f_{c_H\perp c}\right) \\    \Tc_0\left(f_{c_H}\right) & \Tc_0\left(f_{c_H}\right)\\};
  \path[-stealth]
    (m-1-1) edge [double] node [right] {$m$}
    (m-2-1) edge  node [below] {$\overline{f_{c_H\perp c}}$} (m-1-2) 
    (m-2-1.east|-m-2-2) edge node [above] {$\overline{f_{c_H}}$} (m-2-2)
    (m-1-2) edge [double] node [left] {$m$}  (m-2-2);
\end{tikzpicture}
\end{tabular}
\vspace{-0.2cm}
\caption{$\left(\Tc_0\left(f_{c_H}\right),\overline{f_{c_H}}\right)$ is a factor of 
$\left(\Hc\left(f_{c_H}\right),f_{c_H}\right)$ and of $\left(\Tc_0\left(f_{c_H\perp c}\right),\overline{f_{c_H\perp c}}\right)$.}\label{fig:two_diagrams}
\end{figure}
The two diagrams in Figure \ref{fig:two_diagrams} respectively illustrate the two semi-conjugacies $\rho$ and $m$.
With the help of them, one can routinely infer that $\Tc_0\left(f_{c_H}\right)\setminus m(X_H)$ is countable. Therefore, the proof will be completed by the following Lemma \ref{lem:monotone_tuning}.
\end{proof}

\begin{lemm}\label{lem:monotone_tuning}
$m: \Tc_0\left(f_{c_H\perp c}\right)\rightarrow \Tc_0\left(f_{c_H}\right)$
is monotone. Moreover, $m^{-1}(u)$ with $u\in\Tc_0\left(f_{c_H}\right)$  is a singleton if and only if  $\rho^{-1}(u)$ is.
\end{lemm}

\begin{proof}
If $\Tc(f_{c})$ is a point, then $\lambda_{\rm den}\left(f_{c_H\perp c}\right)=\lambda_{\rm den}\left(f_{c_H}\right)$ hence $m$ is a homeomorphism. In the sequel, we assume that $\Tc(f_{c})$ is not a point. 
Observe that   $m$ restricted to $\tau^{c_H\perp c}\left(X_H\right)$ is injective and $\Tc_0\left(f_{c_H}\right)\setminus m\left(\tau^{c_H\perp c}\left(X_H\right)\right)$ is exactly the grand orbit of $\pie(0)$ in $\Big(\Tc_0\left(f_{c_H}\right), \overline{f_{c_H}}\Big)$.
So $m^{-1}(u)$ is a singleton for all  $u\in m\left(\tau^{c_H\perp c}\left(X_H\right)\right)$.

In what follows, we pick a point $u$ in $\Tc_0\left(f_{c_H}\right)\setminus m\left(\tau^{c_H\perp c}\left(X_H\right)\right)$. Then,  $u=\pie(W)$ for some bounded Fatou component $W$. Let $X_u$ denote the point inverse of $u$ under $\tau^{c_H}$. Then $X_u$ is saturated with respect to $\lambda_{\rm den}\left(f_{c_H\perp c}\right)$, in the sense that a class of $\lambda_{\rm den}\left(f_{c_H\perp c}\right)$ is either disjoint from  or entirely contained in $X_u$. There are two observations. First, $X_u$ consists of more than one classes of $\lambda_{\rm den}\left(f_{c_H\perp c}\right)$. Second,   the pre-image of $m^{-1}(u)$ under $\tau^{c_H\perp c}$ equals $A\left(f_{c_H\perp c}\right)\cap X_u$.   By Lemma \ref{lem:convex_hull}, the convex hull ${\rm CH}\Big(\!A\left(f_{c_H\perp c}\right)\Big)$ is the union of all the convex hulls of the classes of  $\lambda_{\rm den}\left(f_{c_H}\right)$ contained in $A\left(f_{c_H\perp c}\right)$. Thus ${\rm CH}\Big(\!A\left(f_{c_H\perp c}\right)\!\Big)\bigcap{\rm CH}(X_u)={\rm CH}\Big(\!A\left(f_{c_H\perp c}\right)\cap X_u\!\Big)$.  Since $m^{-1}(u)=\tau^{c_H\perp c}\Big(A\left(f_{c_H\perp c}\right)\cap X_u\Big)$,  $m^{-1}(u)$ is a continuum by  Lemma \ref{lem:convex_hull_converse}. This proves the monotonicity of $m$.

By the above arguments,  if  $u\in m\left(\tau^{c_H\perp c}\left(X_H\right)\right)$, $m^{-1}(u)$ is a singleton; otherwise, $m^{-1}(u)$ is a non-degenerate continuum and $u=\tau^{c_H}(W)$ for some bounded Fatou component $W$. It follows that $W\subset\rho^{-1}(u)$.  This proves the second part of Lemma \ref{lem:monotone_tuning}. \end{proof}

\begin{rema}\label{rem:monotone_tuning}
In Lemma \ref{lem:monotone_tuning}, the dynamic core $\Big(\Tc_0\left(f_{c_H}\right), \overline{f_{c_H}}\Big)$ of $\left(J_{c_H},f_{c_H}\right)$ is a dendrite factor of $\Big(\Tc_0\left(f_{c_H\perp c}\right), \overline{f_{c_H\perp c}}\Big)$, the dynamic core of $\left(J_{c_H\perp c},f_{c_H\perp c}\right)$. If $\#\Tc\left(f_c\right)=1$, then $m$ is a homeomorphism hence $\Big(\Tc_0\left(f_{c_H}\right), \overline{f_{c_H}}\Big)$ and $\Big(\Tc_0\left(f_{c_H\perp c}\right), \overline{f_{c_H\perp c}}\Big)$  are topologically conjugate. We will use this fact later in proving Theorem \ref{theo:PCF_coincidence}. 
\end{rema}

We define the {\bf central molecule} as follows. Compare \cite{KahnLyubich09-a}.

\begin{defi}\label{def:N_k}
Let $N_0=\overline{H_\heartsuit}$ and  $N_k(k\ge1)$ the union of the closures of all those hyperbolic components $H$ with $\overline{H}\cap N_{k-1}\ne\emptyset$.   Let $M_0$ consist of the atoms of $\M$ that intersect $\overline{\bigcup\limits_{k\ge0}N_k}$ and call it the {\bf central molecule}.
\end{defi}

Applying Theorem \ref{theo:tuning_quadratic} to hyperbolic components $H\subset M_0$, we have the following.

\begin{coro}\label{cor:tuning_molecule}
If $c_H$ is the center  of a hyperbolic component $H\subset M_0$, of period $p\ge2$, then  $h_{\rm GCE}\left(f_{c_H\!\perp\!c}\right)=\frac1p h_{\rm GCE}(f_c)$ holds for all $c\in\M\setminus\{\frac14\}$.	
\end{coro}

By Theorem \ref{theo:tuning_quadratic} and Corollary \ref{cor:tuning_molecule}, we can easily infer the following. 
\begin{coro}\label{cor:tuning_molecule-a}
If $f_c$ with $c\in M_0$ is infinitely renormalizable, then   $h_{\rm GCE}(f_c)=0$. \end{coro}

To conclude this subsection,  we obtain the following. 
\begin{main-theorem}\label{theo:tuning_molecule-b}
A parameter $c\in\M$ belongs to  $\bigcup_{k\ge0}N_k$ if and only if $\#\Tc(f_c)=1$. \end{main-theorem}

Before going on to prove Theorem \ref{theo:tuning_molecule-b}, we recall the basics of orbit portraits  as summarized by Milnor in \cite{Milnor00}.  Given $f\in\Cc_d(d\ge2)$ with filled  Julia set $K$ and   a repelling or parabolic cycle $\Oc=\{z_1,\ldots,z_k\}$,  it is known that there is a finite set $A_i(1\le i\le k)$ such that  the external ray $\Rc(K,\theta)$ lands at $z_i$ if and only if $\theta\in A_i$. Moreover, the map $t\mapsto dt({\rm mod}1)$ acts transitively on the sets $A_1,\ldots, A_k$ and every $\theta\in A_i$ is periodic under $t\mapsto dt$ of the same period, which may not be $k$. By \cite[Theorem 1]{Levin_Przytycki96}, those properties even hold for polynomials having a disconnected Julia set, when the finite set $A_i$ may be replaced by a Cantor set in certain situations.  Following Milnor \cite{Milnor00}, we call $\Pc(\Oc)=\{A_1,\ldots, A_k\}$ the {\bf orbit portrait} of the cycle $\Oc$.

We also refer to \cite{Douady83,Douady87,DH85b,Milnor89} for the basic of the map $c\mapsto c_H\perp c$, which are helpful in describing how the dynamics of $f_c$ vary when $c$ varies over the central molecule $M_0$.

\begin{proof}[{\bf Proof for Theorem \ref{theo:tuning_molecule-b}}]
Let $z_1$ be the  alpha fixed point of  $f_{c_H}$.  Then there exist  $k_0\ge1$ with $k=\frac{p}{k_0}\in\bbZ$ and  $k$ external rays
$\Rc(K_{c_H},\theta_i)(1\le i\le k)$ all of which land at $z_1$.
The forward orbit $\Oc$ of $z_1$ is a singleton $\{z_1\}$ and the orbit portrait $\Pc(\Oc)$ consists of the single set $A_1=\{\theta_1,\ldots,\theta_k\}$.  The union $\{z_1\}\cup \Rc(K_{c_H},\theta_1)\cup\cdots\cup\Rc(K_{c_H},\theta_k)$  cuts the dynamic plane into exactly $k$ open regions. In particular, two of those rays $\Rc(K_{c_H},t_-)$ and $\Rc(K_{c_H},t_+)$, with $t_\pm\in A_1$ and $0<t_-<t_+<1$,  bound the {\bf critical value sector} that contains the critical value $c_H$. Moreover, the two parameter rays $\Rc(\M,t_\pm)$ both land at the root $r'$ of some hyperbolic component $H'$, satisfying $\overline{H'}\cap\overline{H_\heartsuit}=\{r'\}$.  Here $H=H'$ if and only if $k_0=1$.

The union  $\Rc(\M,t_-)\cup\{r'\}\cup \Rc(\M,t_+)$ cuts the parameter plane into two open regions. Denote by $\Wc_{\Pc(\Oc)}$ the one that is disjoint from $H_\heartsuit$.  Note that $c\in\M$ belongs to $\Wc_{\Pc(\Oc)}$  if and only if  $f_c$ has a repelling fixed point with portrait $\Pc(\Oc)$.  See \cite[Theorems 1.1 and 1.2]{Milnor00}.

Let $W_n(0\le n\le p)$ be the Fatou component that contains $v_n=f_{c_H}^n(0)$, with $W_0=W_p$ and $v_0=v_p$. Let $\Theta_n$ consist of all $\theta\in\bbT$ with ${\rm Imp}\left(K_{c_H},\theta\right)\subset\partial W_n$. Each of those sets $\Theta_n$ is a Cantor set. Note that $\sigma_2(\Theta_n)=\Theta_{n+1}$ for $0\le n\le p-1$. In particular, $\sigma_2$ restricts to $\Theta_n$ with $1\le n\le p-1$ is injective
while its restriction to $\Theta_0$ is a two-to-one map. 

Assume in addition that  $c_H\in N_q\setminus N_0(q\ge1)$. Then  the Hubbard tree $\Hc\left(f_{c_H}\right)$ has  no less than $p+1$ vertices, including exactly  $p$ end points $v_0,\ldots, v_{p-1}$.  Therefore, the arcs arc $\overline{v_nz_1}$, connecting $v_n(0\le n\le p-1)$ to $z_1$, are $p$ paths whose union equals the Hubbard tree $\Hc\left(f_{c_H}\right)$.

If $k_0=1$ then $\Hc\left(f_{c_H}\right)$ is a fan having exactly $p$ edges and each  $\overline{v_nz_1}$ is an edge of $\Hc\left(f_{c_H}\right)$. Thus  $\#\Tc\left(f_{c_H}\right)=1$.
If $k_0>1$  the $p$ arcs $\overline{v_nz_1}$ are no longer edges of $\Hc\left(f_{c_H}\right)$,  which then has other branch points than $z_1$.
Either way,  $\#\Tc\left(f_{c_H\perp c}\right)=1$ holds for all $c\ne\frac14$ satisfying $\#\Tc\left(f_c\right)=1$. Since $\#\Tc\left(f_{c}\right)=1$ for all $c\in\overline{H_\heartsuit}$ (See Theorem \ref{theo:major_H}), $\#\Tc\left(f_{c_H\perp c}\right)=1$ still holds for all $c\in\overline{H_\heartsuit}\setminus\{\frac14\}$, hence for all $c\in N_1\setminus\{\frac14\}$. Inductively, one can infer the following.
\begin{lemm}\label{lem:degenerate_core}
Let $c\in M_0$.  If  $f_c$ has a non-repelling cycle then $\#\Tc(f_c)=1$.
\end{lemm}
Recall that for all  $c\ne\frac14$ such that every cycle of $f_c$ is repelling, the alpha fixed point $\alpha_{c_H\perp c}$ of $f_{c_H\perp c}$  cuts $J_{c_H\perp c}$ into exactly $k=\frac{p}{k_0}$ components. It follows that the image $\pie\left(\alpha_{c_H\perp c}\right)$ under $\pie:\widehat{\bbC}\rightarrow\Dc\left(f_{c_H\perp c}\right)$ cuts  $\Tc_0\left(f_{c_H\perp c}\right)$ into exactly $k$ components, whose closures will be denoted by $G_1,\ldots, G_k$.   Notice that each $G_i$ is a dendrite consists of two parts. One is the union of $k_0$ dendrites $\Tc_{i,1},\ldots,\Tc_{i,k_0}\subset\Tc\left(f_{c_H\perp c}\right)$ each of which is isomorphic with  $\Tc_0(f_c)$. The other is the minimal tree  $T_i\subset\Tc\left(f_{c_H\perp c}\right)$ that connects $\pie\left(\alpha_{c_H\perp c}\right)$ to all those $\Tc_{i,j}$.  The following lemma is then immediate, which is the converse of lemma \ref{lem:degenerate_core}.
\begin{lemm}\label{lem:degenerate_core_converse}
Let $c\in M_0$.  If  $\#\Tc(f_c)=1$ then $f_c$ has a non-repelling cycle.
\end{lemm}
Since $c\in M_0$ has a non-repelling cycle if and only if $c$ belongs to $\bigcup_{k\ge0}N_k$. The proof for Theorem \ref{theo:tuning_molecule-b} is completed by Lemmas \ref{lem:degenerate_core} and \ref{lem:degenerate_core_converse}.
\end{proof}

By Corollary \ref{cor:tuning_molecule-a} and Theorem \ref{theo:tuning_molecule-b}, the following is immediate.
\begin{coro}\label{cor:GCE_level_0}
$h_{\rm GCE}(f_c)=0$ for all $c$ belonging to the central molecule $M_0$.	
\end{coro}

\subsection{The GCE Map and the Continuity of its Lower Envelope}\label{section:quadratic_theory}
We will use the growth rate $r_\theta$ for  $\theta\in\bbT=\bbR/\bbZ$ to prove Theorem \ref{theo:properties_of_h}. The continuity of $\log r_\theta$ as a function of $\theta$ is  given in \cite[Theorem 8.2]{Tiozzo16} and independently in \cite{Dudko-Schleicher20}.  We also need the  core decompositions with Peano quotient for planar compacta. 
By \cite[Theorem 7]{LLY19}, every planar compactum $K$ has a core decomposition $\Dc_K^{PC}$, whose elements are called {\bf atoms}. By  definition, if an atom of $K$ intersects the interior $K^o$ it is a singleton. Moreover, $K$ is a Peano compactum  if and only if all its atoms are points. Particularly, if $K$ is a full continuum and $\phi_K:\widehat{\bbC}\setminus K\rightarrow\bbD^*$ the Riemann mapping, then  the prime end impression (shortly, impression) of $K$ at angle 
 $t\in\bbT=\bbR/\bbZ$ coincides with the cluster set of the inverse $\phi^{-1}:\bbD^*\rightarrow\widehat{\bbC}\setminus K$. 
The atoms of $K$ contained in  $\partial K$ give rise to a closed equivalence $\sim_K$ on the unit circle, which identifies $e^{2\pi\textbf{i}s}$ and $e^{2\pi\textbf{i}t}$ if and only if the impressions of $K$ at angles $s,t$ are contained in the same atom of $K$. In particular,  the  equivalence $\sim_M$ coincides with the one given in \cite[Theorem 3]{Douady93}. In such a case, the Riemann mapping $\phi_\M:\widehat{\bbC}\setminus\M\rightarrow\bbD^*$ is often chosen to be the Douady-Hubbard map $\Phi$. Notice that the projection $\pi:\M\rightarrow\Dc_\M^{PC}$ is equivalent to the map $\chi:\M\rightarrow \M_{\rm abs}$, where $\M_{\rm abs}$ is homeomorphic with the pinched disc model given by $\sim_\M$ \cite[IV.6]{Douady93}. In deed, if there is a {\em queer set}, {\em i.e. a point inverse $\chi^{-1}(x)$ that does not reduce to a point}, then $\M$ is not locally connected. If there is a queer set having interior points then every of its interior components are called a {\em queeer component}.  See \cite[Theorem 4]{Douady93}.

Note that the MLC conjecture remains open, although it is known that each atom of $\M$ intersecting $\partial\M$ is the union of finitely many  impressions. Hereafter, for $c\in\M$ let $\delta(c)$ denote the atom of $\M$ that contains $c$ and call $\Rc(\M,t)=\left\{\Phi^{-1}\left(re^{2\pi\textbf{i}t}\right): r>1\right\}$  the external ray at angle $t$. Moreover, we will set $\gamma(t)=\lim\limits_{r\rightarrow\infty}\Phi^{-1}\left(re^{2\pi\textbf{i}t}\right)$, if this limit exists, and  say that the external ray $\Rc(\M,t)$  lands at $\gamma(t)$. We must warn that $\delta(c)$ may properly contain $\{c\}$, even if $\Rc(\M,t)$  lands at $c=\gamma(t)$.

Also note that for any $c_1\ne c_2$  with $\delta(c_1)\ne\delta(c_2)$ one can find a Jordan curve $\Gamma$  such that $\Gamma\cap\M$ is either a point or an arc and the two components of $\M\setminus\Gamma$ respectively contain $c_1,c_2$ \cite{Schleicher04}. On the other hand, if $H$ is a hyperbolic component  with root $r_H\ne\frac14$ then the two external rays landing at $r_H$, together with $r_H$, separates the plane into two domains \cite[Theorem 1.2]{Milnor00}. The one containing $\M_H\setminus\{r_H\}$ is denoted by $\Wc_H$. The other contains the central molecule $M_0$. 
\begin{defi}\label{def:wake}
We call $\Wc_H$ the {\bf wake} at $H$ and $r_H$ the {\bf root point} of $\Wc_H$.	
\end{defi}

The proof for Theorem \ref{theo:properties_of_h} has three steps. First, we compare the values of $h_{\rm GCE}$ and $h_{\rm core}$ at the root $r_H$ and  obtain the following.
\begin{main-theorem}\label{theo:PCF_coincidence}
If  $r_H$ is the root of a hyperbolic component $H$, then $h_{\rm core}(r_H)=h_{\rm GCE}\left(f_{r_H}\right)$.\end{main-theorem}
\begin{proof}
Let $c_H$ be the center of $H$. By \cite[Propositions 18.1 and 18.2]{DH84}, we have $\lambda_{\bbR}\left(f_{r_H}\right)=\lambda_{\bbR}\left(f_{c_H}\right)$ and hence $\lambda_{\rm den}\left(f_{r_H}\right)=\lambda_{\rm den}\left(f_{c_H}\right)$. It follows that $\Big(\Tc_0\left(f_{r_H}\right),\overline{f_{r_H}}\Big)$ is topologically conjugate with $\Big(\Tc_0\left(f_{c_H}\right),\overline{f_{c_H}}\Big)$ and that $h_{\rm GCE}\left(f_{r_H}\right)=h_{\rm GCE}\left(f_{c_H}\right)$. 

By Proposition \ref{cor:tuning_molecule-a} and Theorem \ref{theo:tuning_molecule-b}, we have	$h_{\rm GCE}^{-1}(0)\supset M_0$.  Since every point of $M_0$ may be approximated by an infinite sequence of distinct points in $M_0$, we also have $h_{\rm core}^{-1}(0)\supset M_0$. By Theorem \ref{theo:good_extension}, we have  $h_{\rm GCE}\left(f_{c_H}\right)=h\Big(\Hc\left(f_{c_H}\right),f_{c_H}\Big)$. So, we may assume with no loss of generality that there exist a primitive hyperbolic component $H'$ with center $c_{H'}\notin M_0$ and a hyperbolic component $H''$ with  center  $c_{H''}\in M_0$ that satisfy both  $H\subset\M_{H'}$ and $c_H=c_{H'}\perp c_{H''}$.  
Clearly, the root $r_{H'}$ of $H'$  separates $H$ from $M_0$. 

Let $m: \Tc_0\left(f_{c_{H'}\perp c_{H''}}\right)=\Tc_0\left(f_{c_{H}}\right)\rightarrow \Tc_0\left(f_{c_{H'}}\right)$ be given as in Lemma \ref{lem:monotone_tuning}. Then $m$ is a topological conjugacy between $\Big(\Tc_0\left(f_{c_{H}}\right),\overline{f_{c_{H}}}\Big)$ and $\Big(\Tc_0\left(f_{c_{H'}}\right),\overline{f_{c_{H'}}}\Big)$.  Thus  $h_{\rm GCE}\left(f_{r_{H}}\right)=h_{\rm GCE}\left(f_{c_{H}}\right)=h_{\rm GCE}\left(f_{c_{H'}}\right)=h_{\rm GCE}\left(f_{r_{H'}}\right)$. 

Note that $\M\setminus\{r_{H'}\}$ has exactly two components. The one containing  $c_H$ is exactly $\M\cap\Wc_{H'}$. Moreover, 
$\M\setminus\Big(M_0\cup\{c_{H'}\}\Big)$ has a component $\Uc_{H'}$, satisfying $\overline{\Uc_{H'}}\cap\overline{\Wc_{H'}}=\left\{r_{H'}\right\}$. 

By Theorems \ref{theo:tuning_quadratic} and \ref{theo:tuning_molecule-b},  $h_{\rm GCE}$ is constant on  $\overline{H'}$. See Remark \ref{rem:monotone_tuning}. It follows that $h_{\rm GCE}\left(f_{r_{H'}}\right)$ belongs to the cluster set $C\left(h_{\rm GCE},r_{H'}\right)$.
Since $h_{\rm core}\left(r_{H'}\right)={\rm inf}\ C\left(h_{\rm GCE},r_{H'}\right)$ by definition,  it will suffice to prove that $C\left(h_{\rm GCE},r_{H'}\right)\subset[h_{\rm GCE}\left(f_{r_{H'}}\right),\infty)$. 

 By Theorem \ref{theo:monotonicity}, $h_{\rm GCE}(f_{c})\ge h_{\rm GCE}\left(f_{r_{H'}}\right)$  for all $c$ in 
 $(\M\cap\Wc_{H'})\setminus\M_{H'}$. By Theorem  \ref{theo:tuning_quadratic}, $h_{\rm GCE}\left(f_{c}\right)\ge h_{\rm GCE}\left(f_{c_{H'}}\right)$ for  $c\in\M_{H'}$. It follows that
$h_{\rm GCE}\left(f_{c}\right)\ge h_{\rm GCE}\left(f_{r_{H'}}\right)$ for $c\in \M\cap\overline{\Wc_{H'}}$. 
Therefore, we only need to obtain the following. 

\begin{lemm}\label{lem:coincidence_below}
Given $\{c_n'\}\subset\mathcal{U}_{H'}$ with $c_n'\rightarrow r_{H'}$,  
$h_{\rm GCE}\left(f_{r_{H'}}\right)\le \liminf\limits_{n\rightarrow\infty}h_{\rm GCE}\left(f_{c_n'}\right)$. 
\end{lemm}
 
Let  $t=\liminf\limits_{n\rightarrow\infty}h_{\rm GCE}\left(f_{c_n'}\right)$.
Pick primitive hyperbolic components $H_j(j\ge1)$, with centers $c_j$ and roots $r_j$, that satisfy two requirements: (1)  $c_j\rightarrow r_{H'}$ as $j\rightarrow\infty$, (2) for any $j\ge1$ all but finitely many $c_n'$ are contained in $\Wc_{H_j}$. 
Given $j\ge1$, by Theorem \ref{theo:monotonicity} we can infer that $h_{\rm GCE}\left(f_{c_n'}\right)\ge h_{\rm GCE}\left(f_{c_j}\right)$  for all but finitely many $n$. Since $j\ge1$ is flexible, we can further infer that 
 $t\ge\sup\Big\{h_{\rm GCE}\left(f_{c_j}\right): j\ge1\Big\}.	$
The proof will be completed by the following equation:
\begin{equation}\label{eq:key_equation}
h_{\rm GCE}\left(f_{r_{H'}}\right) = \sup\Big\{h_{\rm GCE}\left(f_{c_j}\right): j\ge1\Big\}.	
 \end{equation}
Pick $\theta_-<\theta_+$ in $\bbT$ such that the external rays $\Rc\left(\M,\theta_\pm\right)$ both land at $r_{H'}$. Similarly, for $j\ge1$ let  $\theta_{j,-}<\theta_{j,+}$ be the two numbers such that  $\Rc\left(\M,\theta_{j,\pm}\right)$ both land at $r_j$. Since $c_j\rightarrow r_{H'}$, we have $\theta_{j,\pm}\rightarrow \theta_\pm$.  Recall that 
 $h_{\rm GCE}\left(f_{c_j}\right)=\log r_{\theta_{j,\pm}}$ for $j\ge1$ and $h_{\rm GCE}\left(f_{r_{H'}}\right)=\log r_{\theta_{\pm}}$. See \cite[Theorem 8.4]{Tiozzo16}. Since $\log r_\theta$ is a continuous function of $\theta$ by \cite[Theorem 8.2]{Tiozzo16}, we shall have  
 $\lim\limits_{j\rightarrow\infty}h_{\rm GCE}\left(f_{c_j}\right)=h_{\rm GCE}\left(f_{r_{H'}}\right)$. This verifies Equation \ref{eq:key_equation} and we are done.
\end{proof}

The second step is to obtain the following.

\begin{main-theorem}\label{theo:monotone_h_core}
The lower envelope $h_{\rm core}:\M\rightarrow\bbR$ is constant in any given component of $\M^o$. It is also constant on any atom of $\M$. Moreover, $h_{\rm core}^{-1}([0,t])$ is continuum for all $t\ge0$.	
\end{main-theorem}
\begin{proof}
The first assertion is immediate, since by definition $h_{\rm GCE}(f_c)$ is constant for $c$ running through any given component of $\M^o$. 

The second one is implied by Theorems \ref{theo:atoms_of_M},  \ref{theo:monotonicity}, \ref{theo:tuning_quadratic} and \ref{theo:PCF_coincidence}. In deed, for any non-degenerate atom $\delta$ of $\M$ there are a sequence of primitive hyperbolic components $H_n$, with roots $r_n$ and periods $p_n\rightarrow\infty$, such that each of the baby Mandelbrot sets $\M_{H_n}$ contains $\delta$. Assume with no loss of generality that $r_n(n\ge1)$ converge to $c_0\in\delta$. Then by Theorems \ref{theo:atoms_of_M}, \ref{theo:tuning_quadratic} and \ref {theo:PCF_coincidence} we have $h_{\rm GCE}(f_c)=\lim\limits_{n\rightarrow\infty}h_{\rm GCE}\left(f_{r_n}\right)=\lim\limits_{n\rightarrow\infty}h_{\rm core}\left(r_n\right)$ for all $c\in\delta$. By Theorem \ref{theo:monotonicity}, we can further verify that the cluster set $C\left(h_{\rm GCE},c\right)$ at any $c\in\delta$ contains $h_{\rm GCE}\left(f_{c_0}\right)$ and is contained in $\Big[h_{\rm GCE}\left(f_{c_0}\right),\log2\Big)$. It then follows that $h_{\rm core}(c)=h_{\rm GCE}\left(f_{c_0}\right)$ for all $c\in\delta$.

The last assertion is to be dealt with as follows. By definition,  $h_{\rm core}$ is lower semi-continuous hence $h_{\rm core}^{-1}([0,t])$ is compact for all $t\ge0$. So we need to verify the connectedness. To do that, we will  employ an argument by contradiction.

Suppose on the contrary that $M_t=h_{\rm core}^{-1}([0,t])$ were disconnected for some $t\ge0$. Then $t<\log2$, since $h_{\rm GCE}(f_c)\le\log2$ for all $c\in\M$ hence $h_{\rm core}^{-1}([0,\log2])=\M$. By Corollary \ref{cor:GCE_level_0}, the central molecule $M_0$ is contained in $h_{\rm core}^{-1}(0)$. Let $P_t$ be the component of $M_t$ containing $M_0$. By the second assertion, we see that if a component  of $M_t$ intersects an atom $\delta$ of $\M$ then it contains $\delta$.  Pick a separation $M_t=E\cup F$ with $P_t\subset E$ and a number $\varepsilon>0$ smaller than the distance between $E$ and $F$. Given $x\in F$, the atom of $\M$ containing $x$ is  contained in $F$. So we can find a hyperbolic component $H_{x}$ whose root $r_x$ separates $x$ from $P_t$ and satisfies 
$\left|x-r_x\right|<\varepsilon$. Due to the choice of $\varepsilon$, we have $r_x\notin E$. By  Theorem \ref{theo:monotonicity}, 
$h_{\rm core}(c)\ge h_{\rm core}\left(r_x\right)$ for all $c\in\Wc_{H_x}$. Thus $h_{\rm core}\left(r_x\right)\le t$ and $r_x\in M_t$. As $\left\{\Wc_{H_{x}}: x\in F\right\}$ is an open cover of $F$, we may pick a finite sub-cover $\Big\{\Wc_{H_{x}}: x\in\{x_1,\ldots,x_k\}\Big\}$.  Then we can find at least one $x\in\{x_1,\ldots,x_k\}$ such that $r_x\notin F$. This is absurd, since $r_x\in M_t=E\cup F$ and $r_{x}\notin E$ by the choices of $\varepsilon$ and  $H_x$. \end{proof}

To prove Theorem \ref{theo:properties_of_h}, the last step is to obtain the following.
\begin{main-theorem}\label{theo:continuity_h_core}
The lower envelope $h_{\rm core}:\M\rightarrow\bbR$ is continuous.	
\end{main-theorem}
\begin{proof}
Suppose on the contrary that  $h_{\rm core}:\M\rightarrow\bbR$ were not uniformly continuous. Then we could find two sequences in $\M$ of distinct points, $\{x_n: n\ge1\}$ and $\{y_n: n\ge1\}$, such that $\lim\limits_{n\rightarrow\infty}x_n=\lim\limits_{n\rightarrow\infty}y_n=c$ and there is a constant $\varepsilon_0>0$ with 
\[\inf\Big\{h_{\rm core}(x_n)-h_{\rm core}(y_n): n\ge1\Big\}\ge\varepsilon_0.
\]
By going to an appropriate subsequence, if necessary, we may assume that $\lim\limits_{n\rightarrow\infty}h_{\rm core}(x_n)=A$ and $\lim\limits_{n\rightarrow\infty}h_{\rm core}(y_n)=B$. The lower semi-continuity of $h_{\rm core}:\M\rightarrow\bbR$ ensures that $B\ge h_{\rm core}(c)$. By Theorems \ref{theo:monotonicity}, \ref{theo:tuning_quadratic}, \ref{theo:PCF_coincidence} and the second assertion of Theorem \ref{theo:monotone_h_core}, we may find a sequence of hyperbolic components $H_j$, whose roots $r_j$ converge to $c$ as $j\rightarrow\infty$, such that $\lim\limits_{j\rightarrow\infty}h_{\rm core}(r_j)= h_{\rm core}(c)$. Similar hyperbolic components may be found for any of the points $x_n$. Thus for any $n\ge1$ we may pick a hyperbolic component $H_n'$, with roots $r_n'$, such that $|r_n'-x_n|$ and $h_{\rm core}(x_n)-h_{\rm core}(r_n')$ are both less than $2^{-n}$. Pick $j_0,n_0\ge1$ such that $h_{\rm core}(x_n)-h_{\rm core}(r_n')<\frac{\varepsilon_0}{2}$ for $n\ge n_0$ and 
$\left|h_{\rm core}(c_j)-h_{\rm core}(c)\right|<\frac{\varepsilon_0}{2}$ for $j\ge j_0$. Then, for any $j\ge j_0$ and any $n\ge n_0$ we have
\begin{equation}\label{eq:h_core_inequality}
	h_{\rm core}(r_n')-h_{\rm core}(r_j)\ge\frac{\varepsilon}{2}.
\end{equation}
Let $\theta_{j,\pm}\in\bbT(j\ge1)$ be the numbers such that the two external rays $\Rc\left(\M,\theta_{j,\pm}\right)$ both land at $r_j$. Let $\Wc_{H_j}$ be the wake at $H_j$, which contains $x$ hence all but finitely many points of the set $\{x_n\}\cup\{y_n\}$. Let $\theta_{n,\pm}'\in\bbT(n\ge1)$ be the numbers such that the two external rays $\Rc\left(\M,\theta_{n,\pm}'\right)$ both land at $r_n'$. Since  $\lim\limits_{n\rightarrow\infty}x_n=\lim\limits_{n\rightarrow\infty}y_n=c$, there is $\theta_0\in\bbT$ such that
\begin{equation}\label{eq:angle_sequences}
\lim\limits_{n\rightarrow\infty}	\theta_{n,\pm}'=\lim\limits_{j\rightarrow\infty}\theta_{j,\pm}=\theta_0.
\end{equation}
By Theorem \ref{theo:PCF_coincidence},  we have $h_{\rm core}(r_j)=\log r_{\theta_{j,\pm}}(j\ge1)$ and $h_{\rm core}(r_n')=\log r_{\theta_{n,\pm}'}(n\ge1)$. It follows that for any $j\ge j_0$ and any $n\ge n_0$ we also have
\begin{equation}\label{eq:r_inequality}
\log r_{\theta_{j,\pm}}	- \log r_{\theta_{n,\pm}'}\ge\frac{\varepsilon}{2}.
\end{equation}
This is absurd, since the map $\theta\mapsto\log r_\theta$ is uniformly continuous \cite[Theorem 8.2]{Tiozzo16}. 
\end{proof}

\begin{rema}\label{rem:further_comments}
By Theorem \ref{theo:PCF_coincidence} and by the constancy of $c\mapsto h_{\rm GCE}(f_c)$ in any interior component of $\M$, $h_{\rm core}(f_c)=h\left(\Hc(f_c),f_c\right)$ for all $c$ such that $f_c(z)$ is PCF. By Theorem \ref{theo:properties_of_h},  $h_{\rm core}:\M\rightarrow\bbR$ continuously extends the classical core entropy map $c\mapsto h(\Hc(f_c),f_c)$, where $c\in\M$ are chosen so that $f_c(z)$ are PCF. Since the set of $c$ with $f_c(z)$ hyperbolic is dense in the real slice $c\in[-2,\frac14]$, by Theorem \ref{theo:good_extension} we also have  $h_{\rm core}(f_c)=h_{\rm GCE}(f_c)=h(\bbR, f_c)$ for $c\in[-2,\frac14]$. 
\end{rema}

\noindent
{\bf Acknowledgement}. The project is supported by (1)  National Key R\&D Program of China (No. 2024YFA1013703); (2) Chinese NSF projects 12571208, 11871483, 12171172, 11901011; (3) Science and Technology Projects of Guangzhou 202102020480; (4) Guangdong Basic and Applied Basic Research Foundation 2021A1515010242; (5) Chinese Postdoctoral Science Foundation No. KLH1411082. The authors want to thank many colleagues for helpful discussions and suggestions, especially Yan Gao, Wen Huang, Wolf Jung, Weiyuan Qiu, Giulio Tiozzo, Yimin Wang, and Jinsong Zeng. In particular, before this project was started, three of the authors joined a mini-course on the  core entropy of polynomials at Peking University by Tiozzo in December, 2019. 

\vspace{1.618cm}

\bibliographystyle{plain}

\end{document}